\documentclass[11pt,leqno]{amsart}

\usepackage[cp1252]{inputenc}
\usepackage[english]{babel}
\usepackage[a4paper]{geometry}
\usepackage{amsmath,amsfonts,amssymb,amsthm,cases,upgreek}
\usepackage{empheq}
\usepackage{stmaryrd}
\usepackage{fullpage}

\usepackage{diagbox}

\usepackage{dsfont}
\usepackage{graphicx}
\usepackage{comment}
\usepackage{bm} 
\usepackage{dsfont}
\usepackage{hyperref}

\newcommand{\lver}{\left\lvert}

\newcommand{\rver}{\right\rvert}

\newcommand\dD{\textrm{d}}

\newcommand\iD{\textrm{i}}

\newcommand{\sgn}{\text{\rm sgn}}

\newcommand\br{\begin{remark}}
\newcommand\er{\end{remark}}
\newcommand\bp{\begin{pmatrix}}
\newcommand\ep{\end{pmatrix}}
\newcommand{\be}{\begin{equation}}
\newcommand{\ee}{\end{equation}}
\newcommand{\bes}{\begin{equation*}}
\newcommand{\ees}{\end{equation*}}
\newcommand\ba{\begin{equation}\begin{aligned}}
\newcommand\ea{\end{aligned}\end{equation}}
\newcommand\bas{\begin{equation*}\begin{aligned}}
\newcommand\eas{\end{aligned}\end{equation*}}

\newcommand{\beg}{\begin{example}}
\newcommand{\eeg}{\end{exaplem}}
\newcommand{\bpr}{\begin{proposition}}
\newcommand{\epr}{\end{proposition}}
\newcommand{\bt}{\begin{theorem}}
\newcommand{\et}{\end{theorem}}
\newcommand{\bc}{\begin{corollary}}
\newcommand{\ec}{\end{corollary}}
\newcommand{\bl}{\begin{lemma}}
\newcommand{\el}{\end{lemma}}
\newcommand{\bd}{\begin{definition}}
\newcommand{\ed}{\end{definition}}
\newcommand{\brs}{\begin{remarks}}
\newcommand{\ers}{\end{remarks}}

\newtheorem{theorem}{Theorem}[section]
\newtheorem{proposition}[theorem]{Proposition}
\newtheorem{corollary}[theorem]{Corollary}
\newtheorem{lemma}[theorem]{Lemma}

\theoremstyle{remark}
\newtheorem{remark}[theorem]{Remark}
\theoremstyle{definition}
\newtheorem{definition}[theorem]{Definition}

\newtheorem{example}[theorem]{Example}

\newcommand\R{\mathbb R}

\newcommand{\N}{\mathbb N}
\newcommand{\Z}{\mathbb Z}

\renewcommand{\footnote}[1]{\textsuperscript{\addtocounter{footnote}{1}(\thefootnote)}\footnotetext{#1}}

\def \epsilon {\varepsilon}

\SetSymbolFont{stmry}{bold}{U}{stmry}{m}{n}

\title{\textbf{Dispersive estimates for nonhomogeneous radial phases: an application to weakly dispersive equations and water wave models}}

\author{Benjamin Melinand}
\address{CEREMADE, CNRS, Universit\'e Paris-Dauphine, Universit\'e PSL, 75016 Paris, France}
\email{{\tt melinand@ceremade.dauphine.fr}}
\thanks{}

\begin{document}

\begin{abstract}
\noindent We study linear dispersive equations in dimension one and two for a class of radial nonhomogeneous phases. $L^{1} \to L^{\infty}$ type estimates, Strichartz estimates, local Kato smoothing and Morawetz type estimates are provided. We then apply our results to different water wave models.
\end{abstract}

\date{\today}
\maketitle

\section{Introduction}

We consider a class of dispersive equations under the form
\be\label{dispersive_eq}
\partial_{t} u = \pm \frac{\iD}{\delta} g(\delta |D|) u
\ee
where $\delta>0$ is a parameter, $u(t,x) \in \mathbb{C}$, $t \in \R$, $x \in \R^{n}$ with $n=1$ or $2$. We assume that

\medskip
(H0) $g$ is a real-valued $\mathcal{C}^{2}$ function defined on $\R^{\ast}_{+}$.
\medskip

In the following we use a Fourier multiplier notation and denote the unique solution of \eqref{dispersive_eq} with $u(0,x)=u_{0}$ well enough localized as 
\be\label{sol_dispersive_eq}
 u(t,x) = e^{\pm \iD \tfrac{t}{\delta} g(\delta |D|)} u_{0} = \frac{1}{(2\pi)^{\frac{n}{2}}} \int_{\R^{n}} e^{\iD x \cdot \xi} e^{\pm \iD \tfrac{t}{\delta} g(\delta |\xi|)} \widehat{u_{0}}(\xi) d\xi.
 \ee
 
\br\label{other_eq_n=1}
When $n=1$, one can also consider dispersive equations under the form
\bes
\partial_{t} u = \pm \frac{g(\delta |D|) }{\delta |D|}\partial_{x} u
\ees
where $\delta > 0$, $u(t,x) \in \R$, $t \in \R$, $x \in \R$ since we have the relation
\[
e^{t \tfrac{g(\delta |D|) }{\delta |D|}\partial_{x}} u_{0} = e^{-\iD \tfrac{t}{\delta} g(\delta |D|)} \mathds{1}_{\{D < 0\}} u_{0} + e^{\iD \tfrac{t}{\delta} g(\delta |D|)} \mathds{1}_{\{D > 0\}} u_{0}.
\]
\er

Studying dispersive estimates for a class of general radial phases is not new and we refer for instance to \cite{decay_estimates_nonhomogeneous_guo_et_al,COX11} or \cite{kenig_ponce_vega_kdv_disper} when $n=1$. We restrict our study to $n=1$ or $n=2$ since we apply our results to water wave models.  Using the methodology developped in this paper together with integration by parts as in \cite{decay_estimates_nonhomogeneous_guo_et_al}, one can easily extend our results to $n \geq 3$, considering $\mathcal{C}^{n}$ phases with ad hoc bounds from above on their derivatives.  Our goal is threefold : (i) provide minimal assumptions on the phase in order to obtain optimal dispersive decays, (ii) study weakly dispersive equations and (iii) obtain bounds that are uniform with respect to $\delta \to 0$. Let us first discuss point (i). The key bounds are of $L^{1} \to L^{\infty}$ type. Using Fourier transform one has to estimate oscillatory integrals. To do so we split the integrals into three different areas: the low, the intermediate and the high frequencies. For the low frequencies, respectively the high frequencies, we consider phases $g$ such that $g''(y) \sim y^{\alpha}$ as $y \to 0$, resp. $y \to \infty$, for some $\alpha \in \R$. The power $\alpha$ can differ between the low and high frequencies. In Lemma \ref{decay1d_1deriv} when $n=1$, Lemma \ref{VdC2d} and Lemma \ref{decay2d_theta=0_HF} when $n=2$, we prove optimal $L^{1} \to L^{\infty}$ type estimates of the low or high frequency part of the propagator. A careful attention is also given for low frequency estimates when $n = 2$ and $g' \gtrsim 1$ at the vicinity of $0$ (Lemma \ref{decay2d_ww_lowfreq}) since a better time decay can be obtained compare to an homogeneous phase. For the intermediate frequencies, the decay depends on the dimension $n$. When $n=1$, $L^{1} \to L^{\infty}$ type estimates with a decay of order $|t|^{-\frac{1}{p}}$ can be obtained if $g'', \cdots , g^{(p)}$ do not vanish at the same time (Lemma \ref{decay1d_multideriv}). When $n = 2$, if $|g'| >0$ and $g'' , \cdots , g^{(p)}$ do not vanish at the same time one can prove $L^{1} \to L^{\infty}$ type estimates with a decay of order $|t|^{-\frac12-\frac{1}{p}}$ (Lemma \ref{decay1d_multideriv}) whereas if $g',g'' , \cdots , g^{(p)}$ do not vanish at the same time one can prove $L^{1} \to L^{\infty}$ type estimates with decay of order $|t|^{-\frac{1}{p}}$ (Remark \ref{decay2d_lowfreq_g'cancels} or Lemma \ref{decay2d_wave_multideriv}). The main tool used to perform these estimates is Van der Corput's Lemma. We also use a Littlewood-Paley decomposition to better understand the regularity needed at low and high frequencies in order to get an optimal time decay. Then, with all these estimates in hand, we perform Strichartz estimates in Subsection \ref{ss:Strichartz}. Finally, in Subsection \ref{ss:kato-morawetz}, we study two properties of dispersive equations: local Kato smoothing effect and local energy decay (related to Morawetz estimates). We provide a result that unifies these two properties.

Concerning point (ii) we have in mind a phase $g$ where $g'(y) \to 0$ as $y \to \infty$ ($\alpha <0$ in the previous paragraph). In that case, extra regularity of the initial datum is necessary in order to obtain a time decay at high frequencies. In Lemma \ref{decay1d_1deriv} when $n=1$, Lemma \ref{VdC2d} and Lemma \ref{decay2d_theta=0_HF} when when $n=2$, we provide the minimal regularity needed to do so. 

Finally to achieve goal (iii) we carefully track the dependence of the parameter $\delta$ in all our dispersive estimates. Assumption (H0) does not necessary mean that $g$ is defined or smooth at $0$. However when it is true with $g'(0) \neq 0$ and $\delta$ small, Equation \eqref{dispersive_eq} can be seen as a perturbation of the half-wave equation. Therefore we expect to capture, at least for the low frequencies, similar properties to the wave equation as $\delta \to 0$. We have in mind $L^{1} \to L^{\infty}$ type estimates with decay of order $|t|^{-\frac{1}{2}}$ when $n=2$ and local energy decay for any dimension (also called Morawetz estimates \cite{Morawetz68}). We show such properties in Lemma \ref{decay2d_wave_2deriv} (for the $L^{1} \to L^{\infty}$ type estimates) and in Proposition \ref{kato-morawetz} and Corollary \ref{decay_local_energy} (for the bound and the decay of local energy) with bounds that are uniform with respect to $\delta \to 0$.

This problem is motivated by the study of water-wave models in the shallow water regime. Adopting the framework of \cite{Lannes_ww}, we introduce the shallowness parameter $\mu \in (0,1]$. As examples of our dispersive estimates, we mainly study two models. Firstly, in Subsection \ref{ss:ww}, we consider the linearized irrotational water wave equations (under the Zakharov/Craig-Sulem formulation)
 \begin{equation}\label{linear_ww}\left\{\begin{array}{l}
\partial_t\zeta - \frac{\tanh(\sqrt{\mu}|D|)}{\sqrt{\mu}} |D| \psi=0,\\
\partial_t\psi+\zeta=0,
\end{array}\right.\end{equation}
so that $\delta=\sqrt{\mu} \in (0,1]$ and 
\[
g(y) = \sqrt{y \tanh(y)}.
\]
Secondly, we study in Subsection \ref{ss:abcd_Boussi} the dispersive part of the linearized abcd-Boussinesq systems 
 \begin{equation}\label{linear_abcd}\left\{\begin{array}{l}
(1-\mu b \Delta)\partial_t\zeta + (1+\mu a \Delta) \nabla \cdot V =0,\\
(1-\mu d \Delta)\partial_t \nabla \cdot V+(1+c \Delta) \nabla \zeta=0,
\end{array}\right.\end{equation}
under the condition that $b \geq 0$, $d \geq 0$, $a \leq 0$, $c \leq 0$ (in order to get the wellposedness \cite{bona_chen_saut_derivation}), so that  $\delta=\sqrt{\mu} \in (0,1]$ and 
\[
g(y) = y \sqrt{\frac{(1-\mu a y^2)(1-\mu c y^2)}{(1+\mu b y^2)(1+\mu d y^2)}}.
\]
These two models exactly fall within our goals (ii)  (at least for most of parameters $a,b,c,d$ when one considered System \eqref{linear_abcd}) and (iii) (since $\mu$ can be very small). 

We also provide dispersive estimates for other models when $n=1$:

\begin{itemize}
\item[$-$] the linearized Ostrovsky equation, studied in Subsection \ref{ss:ostrovsky},
\[
\partial_{t} u = ( - \partial_{x}^{-1} u + b \partial_{x}^{3} u )
\]
where, in the setting of Remark \ref{other_eq_n=1}, $\delta=1$ and $g(y) = \frac{1}{y} - b y^{3}$,

\item[$-$] a linearized BBM/KdV equation, studied in Subsection \ref{ss:BBM_KdV},
\[
\partial_{t} u + \mu p \partial_{x}^{2} \partial_{t} u = \pm (\partial_{x} u + \mu (p+\tfrac16) \partial_{x}^{3} u)
\]
with $p \leq 0$ and where, in the setting of Remark \ref{other_eq_n=1}, $\delta=\sqrt{\mu}$ and $g(y) =y \frac{1 -  (p+\tfrac16) y^2}{1- p y^2}$,

\item[$-$] the linear intermediate long wave equation, studied in Subsection \ref{ss:ILW}, 
\[
u_{t} =  \frac{1}{\rho} \varphi(\rho |D|) \partial_{x} u
\]
where $\varphi(y) = y \coth(y) - 1$ and with $\rho>0$. 
\end{itemize}

\subsection*{Notations}
\begin{itemize}
\item If $u$ is a Schwartz class function, we define $\mathcal{F}u$ or $\hat{u}$ as the Fourier transform of $u$ by
\[
\hat{u}(\xi) = \frac{1}{(2\pi)^{\frac{n}{2}}} \int_{\R^{n}} e^{-\iD x \cdot \xi} u(x) dx.
\]

\item If $f$ is a smooth function that is at most polynomial at infinity, we define the Fourier multiplier $f(D)$ as, for any Schwartz class function $u$, 
\[
f(D) u = \mathcal{F}^{-1} (  f(\xi) \hat{u}(\xi)).
\]

\item The family $(Q_{j})_{j \in \Z}$ is defined in Subsection \ref{ss:LittlewoodPaley}.

\item $J_{0}$ is a Bessel function defined as $J_{0}(s) = \int_{0}^{2\pi} e^{is \sin(\theta)} d\theta$.

\item If $g$ is a function and $p \in \mathbb{N}$, we denote by $g^{(p)}$ the $p$-th derivative of $g$.

\item $\Z_{J} := \{ k \in \mathbb{Z} \text{  , } [2^{k-1},2^{k+1}] \subset J \}$ where $J$ is an interval.

\item If $p \in [1,\infty]$, we denote $p' = \tfrac{p}{p-1}$.

\item If $A$ is  a subset of $\R^{n}$, the map $z \in \R^{n}  \to \mathds{1}_{A}(z)$ stands for the indicator function of $A$.

\item If $u_{0}$ is function defined on $\R^{n}$ we denote by $\| u_{0} \|_{L^{p}}$ the $L^{p}(\R^{n})$ norm of $u_{0}$.

\item If $F : (t,x) \in \R \times \R^{n} \to F(t,x) \in \R$, the norm $\| F \|_{L^{q}_{t} L^{r}_{x}}$ corresponds to the norm of the space $L^{q}(\R;L^{r}(\R^{n}))$. 

\item If $T : E \to F$ is linear bounded operator with $E,F$ two Banach spaces, we denote by $T^{\ast}$ its adjoint.
\end{itemize}

\subsection*{Acknowledgments} We thank Jean-Claude Saut for his interest in this work and Vincent Duch\^ene for his comments that greatly improve the presentation of this paper. We also thank Anatole Gaudin for very enriching discussions about interpolation theory and harmonic analysis. This work has been partially funded by the ANR project CRISIS (ANR-20-CE40-0020-01).

\section{Dispersive estimates}\label{sec_dispersive}

 \subsection{Littlewood-Paley decomposition}\label{ss:LittlewoodPaley}
Since the phases we consider are not necessary homogeneous or weakly dispersive, we use a Littlewood-Paley decomposition to better capture the difference between the low and the high frequencies and catch possible need of extra regularity. We introduce a smooth nonnegative even function $\varphi_{0}$ supported in $[-\tfrac45,\tfrac45]$, that is equal to $1$ in $[-\tfrac35,\tfrac35]$ and that is nonincreasing on $\R^{+}$. Then we define, for any $y \in \R$ and any $j \in \mathbb{Z}$, $Q_{j}(y) :=\varphi_{0}(2^{-j-1} y)-\varphi_{0}(2^{-j} y)$. We note that $Q_{j}$ is a function supported in the annulus $\mathcal{C}( \tfrac35 2^{j}, \tfrac85 2^{j}) \subset \mathcal{C}( 2^{j-1}, 2^{j+1})$ and that for any $y \in \R^{\ast}$
\[
Q_{j}(y) \in [0,1] \text{   ,   } \sum_{j \in \mathbb{Z}} Q_{j}(y) = 1 \text{   ,   } \frac12 \leq \sum_{j \in \mathbb{Z}} Q^{2}_{j}(y) \leq 1.
\]

We also define the set $\Z_{J}$ for any interval $J$ of $\R^{+}$ as
\[
\Z_{J} := \{ k \in \mathbb{Z} \text{  , } [2^{k-1},2^{k+1}] \subset J \}.
\]
Roughly speaking for some $y_{1}>y_{0}>0$, $\Z_{J}$ gathers the low frequencies when $J=(0,y_{0}]$, the intermediate frequencies when $J=[y_{0},y_{1}]$ and the high frequencies when $J=[y_{1},\infty)$. Note that if $J=(y_{0},y_{1})$ with $0 \leq y_{0} \leq \tfrac{3}{32} y_{1} \leq \infty$, then $\sum_{j \in \Z_{J}} Q_{j}$ is equal to $1$ on $(\tfrac{16}{5} y_{0}, \tfrac{3}{10} y_{1})$.

\subsection{Tools for dispersive estimates and framework}\label{tools_disp}

In this section, we introduce some basic tools in order to prove decay estimates. We first recall the van der Corput lemma (see for instance \cite{Stein_harmonic}). 
\begin{lemma}[Van der Corput]
Let $a<b$ be real numbers, $\lambda>0$ and $\phi$ a smooth real-valued function defined on $(a,b)$. Assume that there exists $p \in \N^{\ast}$ such that $|\phi^{(p)}| \geq \lambda$ on $(a,b)$ and, if $p=1$, that $\phi'$ is monotonic on $(a,b)$. Then, there exists a constant $c_{p}$ that only depends on $p$, such that for any $t \in \R^{\ast}$ and any $\mathcal{C}^1$ function $f$ defined on $[a,b]$
\[ \lver \int_{a}^{b} e^{i t \phi(x)} f(x) dx \rver \leq c_{p} (\lambda |t|)^{-\frac{1}{p}} \left( |f(b)| + \int_{a}^{b} |f'(x)| dx\right).\]
\end{lemma}
We easily get the following corollary in the case $p=1$.
\begin{corollary}\label{cor_vandercorput}
Let $a<b$ be real numbers, $r \in \N$, $\lambda>0$ and $\phi$ a smooth real-valued function defined on $(a,b)$. Assume that $|\phi'| \geq \lambda$ on $(a,b)$ and that $\phi''$ has at most $r$ zeros on $(a,b)$. Then, there exists a constant $c$ such that for any $t \in \R^{\ast}$ and any smooth function $f$ defined on $[a,b]$
\[ \lver \int_{a}^{b} e^{i t \phi(x)} f(x) dx \rver \leq \frac{c}{\lambda |t|} \left( r \| f \|_{\infty} + |f(b)| + \int_{a}^{b} |f'(x)| dx\right).\]
\end{corollary}

Then, we introduce the Bessel function $J_{0}$: for any $s \in \R$
\[
J_{0}(s) = \int_{0}^{2\pi} e^{is \sin(\theta)} d\theta.
\]
A direct application of Van der Corput's Lemma gives
\[
 |J_{0}(s)| + |J_{0}'(s)| \lesssim \frac{1}{\sqrt{1+|s|}}.
\]
Furthermore, introducing a smooth nonnegative function $\chi$ defined on $[0,2\pi] $ with $\chi(x)=0$ for $x\in [\frac{5\pi}{4}, \frac{7 \pi}{4}]$ and $\chi(x)=1$ for any $ x\in [\frac{\pi}{4}, \frac{3 \pi}{4}]$, one can decompose $J_{0}$ as
\[
J_{0}(s) = e^{is} \int_{0}^{2\pi} \chi(\theta) e^{is (\sin(\theta)-1)} d\theta + e^{-is} \int_{0}^{2\pi} (1-\chi(\theta)) e^{is (\sin(\theta)+1)} d\theta := e^{is} h_{-}(s) + e^{-is} h_{+}(s).
\]
Integrating by parts if necessary and using Van der Corput's Lemma, one can get for any $p \in \N$
\[
|h_{\pm}^{(p)}(s)| \lesssim \frac{1}{(1+|s|)^{p+\frac12}}.
\]

In the next two subsections, we prove different decay estimates. Our goal is to provide $L^{1} \to L^{\infty}$ type bounds of \eqref{sol_dispersive_eq} and therefore $L^{\infty}_{x}$ bounds of 
\[
\int_{\R^{n}} e^{\iD x \cdot \xi} e^{\iD \tfrac{t}{\delta} g(\delta |\xi|)} |\xi|^{s} d\xi.
\]
The quantity $|\xi|^{s}$ is introduced to catch global smoothing effects or to assume extra regularity in order to get a time decay. Since $g$ can have different behaviors on the low, intermediate and high frequencies, we split the domain of integration into different pieces and it remains to estimate two type of integrals:
\[
\int_{\R^{n}} e^{\iD x \cdot \xi} e^{\iD \tfrac{t}{\delta} g(\delta |\xi|)} \chi( \delta |\xi|) |\xi|^{s} d\xi
\]
where $\chi$ is a smooth function supported on an interval and
\[
\int_{\R^{n}} e^{\iD x \cdot \xi} e^{\iD \tfrac{t}{\delta} g(\delta |\xi|)} P( \tfrac{\delta |\xi|}{2^{j}}) |\xi|^{s} d\xi
\]
for some $j \in \Z$ and where $P$ is a smooth function whose support is contained in the annulus $\mathcal{C}(\tfrac12,2)$. The first type of integrals appears for the intermediate frequencies or when a derivative of $g'$ is far from $0$. The second type of integrals occurs when one performs a Littlewood-Paley decomposition and wants to understand the regularity needed to get an optimal time decay.

\subsection{Dispersive estimates in the 1d case}\label{ss:n=1}

In this subsection we assume that $n=1$. We first give an easy consequence of Van der Corput's Lemma in case a derivative of $g$ is far from $0$ on an interval.

\bl\label{decay1d_VdC}
Let $J \subset \R^{+}$ an interval. Let $\lambda> 0$, $p \in \mathbb{N}$ with $p \geq 2$. Assume that $g$ satisfies (H0), that $g$ is $\mathcal{C}^{p}$ and $|g^{(p)}| \geq \lambda$ on $J$. Let $\chi$ be a smooth function on $\R$ whose support is a subset of $J$ and such that $\chi' \in L^{1}(J)$. There exists $C>0$ such that for any $t \in \R^{\ast}$ and any $\delta>0$
\[
\sup_{x \in \R} \left| \int_{\R} e^{\iD x \xi} e^{\iD \tfrac{t}{\delta} g(\delta |\xi|)} \chi(\delta \xi) d\xi \right| \leq C \frac{\delta^{\frac{1-p}{p}}}{|t|^{\frac{1}{p}}}.
\]
\el

\br 
When  $J$ is unbounded, the previous integral is well defined from integration by parts since $\frac{g''(y)}{g'(y)^2}$ has a sign for $|y|$ large enough and $\displaystyle \lim_{|y| \to +\infty} |g'(y)| = + \infty$.
\er

The previous lemma is particularly useful for homogeneous phases or low frequency estimates. However, in many situations, one has to be more accurate especially when $g'(y)\to 0$ as $y \to +\infty$. For this reason we localize in the annulus $\mathcal{C}(\delta 2^{k-1}, \delta 2^{k+1})$ for some $k \in \Z$. Fix $P$ a smooth even function supported in the annulus $\mathcal{C}(\tfrac12,2)$. We define the quantity for any $t \in \R$, $x \in \R$, $s \in \R$ and $k \in \mathbb{Z}$
\[
I^{s}_{t,x,k} := \int_{\R^{+}} e^{\iD x \xi} e^{\iD \tfrac{t}{\delta} g(\delta |\xi|)} P(\tfrac{\delta \xi}{2^{k}}) \xi^{s} d\xi = \frac{2^{(1+s)k}}{\delta^{s+1}} \int_{\frac12}^{2} e^{\iD \frac{2^{k}}{\delta} x r} e^{ \iD \tfrac{t}{\delta} g(2^{k} r)} P(r) r^{s} dr.
\]
Note that one can easily compute from $I^{s}_{t,x,k}$ the quantities
\[
\int_{\R} e^{\iD x \xi} e^{\iD \tfrac{t}{\delta} g(\delta |\xi|)} P(\tfrac{\delta \xi}{2^{k}}) |\xi|^{s} d\xi \text{  or  } \int_{\R} e^{\iD x \xi} e^{\iD t \tfrac{g(\delta |\xi|)}{\delta |\xi|} \xi} P(\tfrac{\delta \xi}{2^{k}})\xi^{\alpha} d\xi \text{   for } \alpha \in \N.
\]

\bl\label{decay1d_1deriv}
Let $J \subset \R^{+}$ an interval. Let $\lambda > 0$, $\alpha \in \R$, $s \in \R$ and $p \in \N$ with $p \geq 2$. Assume that $g$ satisfies (H0), that $g$ is $\mathcal{C}^{p}$ and that $|g^{(p)}(y)| \geq \lambda y^{\alpha}$ on $J$. There exists $C>0$ such that for any $t \in \R^{\ast}$, any $\delta>0$ and any $k \in \Z_{J}$
\[
\sup_{x \in \R} |I^{s}_{t,x,k}| \leq C \frac{2^{(s-\frac{\alpha}{p}) k} \delta^{\frac{1-sp-p}{p}}}{|t|^{\frac{1}{p}}},
\]
and, if $l \geq 2$ and $|g''(y)| \geq \lambda y^{\alpha}$ on $J$,
\[
\sup_{x \in \R} |I^{s}_{t,x,k}| \leq C \frac{2^{(s+\frac{l-2-\alpha}{l}) k} \delta^{\frac{1-sl-l}{l}}}{|t|^{\frac{1}{l}}}.
\]
\el

\begin{proof}
By a change of variables, one can assume that $\delta=1$. Defining $\phi(r)= t g(2^{k} r) + 2^{k} x r$, we have $|\phi^{(p)}(r)| = 2^{pk} |t g^{(p)}(2^{k} r)|$. Van der Corput's Lemma gives
\[
|I^{s}_{t,x,k}| \lesssim \frac{2^{(1+s)k}}{(2^{pk} |t| g^{(p)}(2^{k} r_{1}))^{\frac{1}{p}}}
\]
for some $r_{1} \in [\tfrac12,2]$. The first bound follows. The second bound can be obtained by interpolation between the first bound and the fact that $|I^{s}_{t,x,k}| \lesssim 2^{(s+1)k}$.
\end{proof}

By taking $P=Q_{0}$ in the previous lemma we get
\[
\sup_{x \in \R} \left| \int_{\R} e^{\iD x \xi} e^{\iD \tfrac{t}{\delta} g(\delta |\xi|)} \sum_{k \in \Z_{J}} Q_{j}(\delta \xi) |\xi|^{s} d\xi \right| \lesssim \frac{\delta^{\frac{1-sp-p}{p}}}{|t|^{\frac{1}{p}}}
\]
provided that $s-\frac{\alpha}{p}>0$ and $J=(0,y_{0})$ or $s-\frac{\alpha}{p}<0$ and $J=(y_{0},+\infty)$. 

The next lemma can be seen as an alternative to Lemma \ref{decay1d_VdC} that is better suited for inhomogeneous phases when $g''(y) \to 0$ as $y \to 0$ or $+\infty$ without vanishing.

\bl\label{decay1d_VdC2}
Let $J \subset \R^{+}$ an interval. Let $\Lambda,\lambda> 0$, $\ell ,s, \alpha \in \R$ with $\alpha \neq -2$. Assume that $g$ satisfies (H0) and that $|g''(y)| \geq \lambda y^{\alpha}$ on $J$. 
\begin{enumerate}
\item If $\frac{\alpha}{2} > s > -1$  or $\frac{\alpha}{2} <s < -1$, there exists $C>0$ such that for any $t \in \R^{\ast}$ and any $\delta>0$
\[
\sum_{k \in \Z_{J}} \sup_{x \in \R} \left| I^{s}_{t,x,k} \right| \leq C \frac{\delta^{-\frac{(1+\alpha)(s+1)}{2+\alpha}}}{|t|^{\frac{s+1}{2+\alpha}}} \text{     ,     } \sup_{x \in \R} \left| \int_{\R} e^{\iD x \xi} e^{\iD \tfrac{t}{\delta} g(\delta |\xi|)} \sum_{k \in \Z_{J}} Q_{k}(\delta \xi)  |\xi|^{s} d\xi \right| \leq C \frac{\delta^{-\frac{(1+\alpha)(s+1)}{2+\alpha}}}{|t|^{\frac{s+1}{2+\alpha}}}.
\]
\item If furthermore $\Lambda y^{\alpha+1} \geq |g'(y)-\ell| \geq \lambda y^{\alpha+1}$ on $J$ and $\alpha \notin \{-2,-1\}$, there exists $C>0$ such that for any $t \in \R^{\ast}$ and any $\delta>0$
\[
\sup_{x \in \R} \;\; \sum_{k \in \Z_{J}} \left| I^{\frac{\alpha}{2}}_{t,x,k} \right| \leq C \frac{\delta^{-\frac{1+\alpha}{2}}}{\sqrt{|t|}} \text{     ,     } \sup_{x \in \R} \left| \int_{\R} e^{\iD x \xi} e^{\iD \tfrac{t}{\delta} g(\delta |\xi|)} \sum_{k \in \Z_{J}} Q_{k}(\delta \xi) |\xi|^{\frac{\alpha}{2}} d\xi \right| \leq C  \frac{\delta^{-\frac{1+\alpha}{2}}}{\sqrt{|t|}}.
\]
\end{enumerate}
\el

\br 
When $J=]0,1]$, $\delta=1$ and $s=0$, we get item (c) of Theorem 1 in \cite{decay_estimates_nonhomogeneous_guo_et_al}. Point (2) is similar to Corollary 2.9 in \cite{kenig_ponce_vega_kdv_disper}. 

When $J=(0,y_{0}]$ and $\frac{\alpha}{2} <s < -1$ or $J=[y_{0},\infty)$ and $\frac{\alpha}{2} > s > -1$ for some $y_{0}>0$, we rather use Lemma \ref{decay1d_1deriv} since the decay provided by Lemma \ref{decay1d_VdC2} is not optimal.
\er

\br\label{r:alphanot2}
As noted in \cite{kenig_ponce_vega_kdv_disper}, the assumption $\alpha \neq -2$ is crucial since the estimate is not true for $g(\xi)=\ln(|\xi|)$. Concerning the case $(\alpha,s)=(-1,-\frac12)$, our proof fails in the region $|x| \sim |t|$ since the set $A_{t,x}$ defined in the proof can actually be $\Z$ and then is not bounded. We interpret this issue as the fact that $\alpha=-1$ corresponds to a wave-type behavior so that the dispersive effects can not be optimal.
\er

\begin{proof}
We begin with the first estimate. In that case, $(s+1)(s-\frac{\alpha}{2})<0$. We have, using Lemma \ref{decay1d_1deriv}, for any $k \in \Z_{J}$
\[
|I^{s}_{t,x,k}| \lesssim \min \left( 2^{(1+s)k} \delta^{-(s+1)} , \frac{2^{(s-\frac{\alpha}{2}) k} \delta^{-(s+\frac12)}}{\sqrt{|t|}} \right) 
\]
so that if $\alpha>-2$
\[
\sum_{k \in \Z_{J}} \sup_{x \in \R} \left| I^{s}_{t,x,k} \right| \lesssim \sum_{2^{k} |t|^{\frac{1}{\alpha+2}} \leq \delta^{\frac{1}{\alpha+2}}} 2^{(1+s)k} \delta^{-(s+1)}  + \sum_{2^{k} |t|^{\frac{1}{\alpha+2}} \geq \delta^{\frac{1}{\alpha+2}}} \frac{2^{(s-\frac{\alpha}{2}) k} \delta^{-(s+\frac12)}}{\sqrt{|t|}},
\]
whereas if $\alpha<-2$
\[
\sum_{k \in \Z_{J}} \sup_{x \in \R} \left| I^{s}_{t,x,k} \right| \lesssim \sum_{2^{k} |t|^{\frac{1}{\alpha+2}} \leq \delta^{\frac{1}{\alpha+2}}} \frac{2^{(s-\frac{\alpha}{2}) k} \delta^{-(s+\frac12)}}{\sqrt{|t|}} + \sum_{2^{k} |t|^{\frac{1}{\alpha+2}} \geq \delta^{\frac{1}{\alpha+2}}} 2^{(1+s)k} \delta^{-(s+1)} .
\]
The first estimate follows and also the second one by taking $P=Q_{0}$ in the definition of $I^{s}_{t,x,k}$. Then, we consider the third estimate. By translation, one can assume that $\ell=0$. Two cases occur.

\medskip
\noindent$\bullet$ Case 1: $\frac{\lambda}{4} 2^{(\alpha+1) k} |t| \leq |x| \leq 4 \Lambda 2^{(\alpha+1) k} |t|$
\medskip

Lemma \ref{decay1d_1deriv} gives
\[
|I^{\frac{\alpha}{2}}_{t,x,k}| \lesssim \frac{\delta^{-\frac{1+\alpha}{2}}}{\sqrt{|t|}}.
\]
Note that since $\alpha \neq -1$, the set $A_{t,x} := \{ k \in \mathbb{Z}_{J}, \frac{\lambda}{4} 2^{(\alpha+1) k} |t| \leq |x| \leq 4 \Lambda 2^{(\alpha+1) k} |t| \}$ is bounded by a number independent of $t,x,\delta$.

\medskip
\noindent $\bullet $ Case 2: $\frac{\lambda}{4} 2^{(\alpha+1)k} |t| \geq |x|$ or $|x| \geq 4 \Lambda 2^{(\alpha+1) k} |t|$
\medskip

Defining the phase $\phi(r)= \tfrac{t}{\delta} g(2^{k}r) + \frac{2^{k}}{\delta} x r$ and since $|\phi'(r)| \geq \frac{\lambda}{4} \delta^{-1} 2^{(\alpha+2)k} |t|$, we get thanks to Corollary \ref{cor_vandercorput} 
\[
|I^{\frac{\alpha}{2}}_{t,x,k}| \lesssim \min \left( 2^{\frac{2+\alpha}{2} k} \delta^{-\frac{2+\alpha}{2}} , \frac{2^{-\frac{\alpha+2}{2} k} \delta^{-\frac{\alpha}{2}}}{|t|} \right).
\]
Therefore gathering all previous estimates, we get if $\alpha>-2$ and $\alpha \neq -1$
\begin{align*}
\sum_{k \in \Z_{J}} \left|  I^{\frac{\alpha}{2}}_{t,x,k} \right| &\lesssim \sum_{k \in A_{t,x}} | I^{\frac{\alpha}{2}}_{t,x,k}| + \sum_{k \in \Z_{J} \text{, } k \notin A_{t,k}} | I^{\frac{\alpha}{2}}_{t,x,k}|\\
&\lesssim \frac{\delta^{-\frac{1+\alpha}{2}}}{\sqrt{|t|}} + \sum_{2^{k} |t|^{\frac{1}{\alpha+2}} \leq \delta^{\frac{1}{\alpha+2}}} 2^{\frac{2+\alpha}{2} k} \delta^{-\frac{2+\alpha}{2}}  + \sum_{2^{k} |t|^{\frac{1}{\alpha+2}} \geq \delta^{\frac{1}{\alpha+2}}} \frac{2^{-\frac{\alpha+2}{2} k} \delta^{-\frac{\alpha}{2}}}{|t|},
\end{align*}
whereas if $\alpha <-2$
\[
\sum_{k \in \Z_{J}} \left|  I^{\frac{\alpha}{2}}_{t,x,k} \right| \lesssim \frac{\delta^{-\frac{1+\alpha}{2}}}{\sqrt{|t|}} + \sum_{2^{k} |t|^{\frac{1}{\alpha+2}} \leq \delta^{\frac{1}{\alpha+2}}} \frac{2^{-\frac{\alpha+2}{2} k} \delta^{-\frac{\alpha}{2}}}{|t|} + \sum_{2^{k} |t|^{\frac{1}{\alpha+2}} \geq \delta^{\frac{1}{\alpha+2}}} 2^{\frac{2+\alpha}{2} k} \delta^{-\frac{2+\alpha}{2}} .
\]
The third and fourth estimates follow.

\end{proof}

Finally, we provide a result for the intermediate frequencies. We consider the situation where some derivatives of $g$ do not vanish at the same time.

\bl\label{decay1d_multideriv}
Let $\lambda > 0$, $y_{1} \geq y_{0} > 0$ and $l \in \N$ with $l \geq 2$. Assume that $g$ is $\mathcal{C}^{l}$ and that $\displaystyle \sum_{p=2}^{l} |g^{(p)}|\geq \lambda$ on $[\tfrac12 y_{0}, 2 y_{1}]$. There exists $C>0$ such that for any $k \in \mathbb{Z}$ such that $y_{0} \leq 2^{k} \leq y_{1}$, any $\delta>0$ and any $t \in \R^{\ast}$
\[
\sup_{x \in \R} |I^{s}_{t,x,k}| \leq C \frac{\delta^{-(s+1-\frac{1}{l})k}}{|t|^{\frac{1}{l}}} \text{        ,        } \sup_{x \in \R} \left| \int_{\R} e^{\iD x \xi} e^{\iD \tfrac{t}{\delta} g(\delta |\xi|)} \sum_{y_{0} \leq 2^{k} \leq y_{1}} Q_{k}(\delta \xi)  |\xi|^{s} d\xi \right| \leq C \frac{\delta^{-(s+1-\frac{1}{l})k}}{|t|^{\frac{1}{l}}}.
\]
\el

\begin{proof}
For $p \in \{2,\cdots,l \}$ we define the sets
\[
J_{p} := \left\{y \in [\tfrac12 y_{0}, 2 y_{1}], |g^{(p)}(y)| > \frac{\lambda}{2(l-1)} \right\} \bigcap_{k=2}^{p-1} \left\{ y \in [\tfrac12 y_{0}, 2 y_{1}], |g^{(k)}(y)| \leq \frac{\lambda}{2(l-1)} \right\}
\]
so that by assumption $J_{p}$ is a finite union of intervals and $\sqcup_{p=2}^{l} J_{p} = [\tfrac12 y_{0}, 2 y_{1}]$, and the integrals
\[
I_{p} :=  \frac{2^{(s+1)k}}{\delta^{s+1}} \int_{\frac12}^{2} e^{ \iD \frac{2^{k}}{\delta} x r} e^{ \iD \tfrac{t}{\delta} g(2^{k} r)} \mathds{1}_{J_{p}}(2^{k} r) P(r) r^{s} dr.
\]
If $2^{k} |t| \leq \delta$
\[
|I_{p}| \lesssim \frac{2^{(s+1)k}}{\delta^{s+1}} = \frac{2^{(s+1-\frac{1}{l})k}}{\delta^{(s+1-\frac{1}{l})k}} \frac{2^{\frac{k}{l}}}{\delta^{\frac{k}{l}}} \leq \frac{2^{(s+1-\frac{1}{l})k} \delta^{-(s+1-\frac{1}{l})k}}{|t|^{\frac{1}{l}}}.
\]
If now $2^{k} |t| \geq \delta$, denoting $\phi(r) = \tfrac{t}{\delta} g(2^{k} r) + \frac{2^{k}}{\delta} xr$, $|\phi^{(p)}(r)| = 2^{p k} \delta^{-1} |t g^{(p)}(2^{k} r)|$ and by Van der Corput's Lemma
\[
|I_{p}| \lesssim \frac{2^{sk}}{\delta^{s+1}} \left( \frac{\delta}{|t|} \right)^{\frac{1}{p}} \lesssim 2^{sk} \frac{\delta^{-s-\frac12}}{\sqrt{|t|}} + 2^{sk} \frac{\delta^{-(s+1-\frac{1}{l})}}{|t|^{\frac{1}{l}}}.
\]
Then since $2^{k} |t| \geq \delta$
\[
\frac{\delta^{-s-\frac12}}{\sqrt{|t|}}  = \frac{\delta^{-s-\frac12}}{|t|^{\frac{1}{l}}} \frac{1}{|t|^{(\frac{1}{2}-\frac{1}{l})}} \leq \frac{\delta^{-(s+1-\frac{1}{l})} 2^{(\frac{1}{2}-\frac{1}{l})k}}{|t|^{\frac{1}{l}}}.
\]
The first bound follows noticing that $2^{k} \sim 1$. The second bound follows by taking $P=Q_{0}$ in the definition of $I^{s}_{t,x,k}$ and the boundedness of the domain of summation.
\end{proof}

\subsection{Dispersive estimates in the 2d case}\label{ss:n=2}

We assume that $n=2$. We begin with a series of results that provide a better decay compare to the wave equation.  There are however not uniform with respect to $\delta \to 0$. First, we consider the quantity for $t \in \R$, $x \in \R^{2}$ and $\chi$ a smooth function
\begin{equation}\label{I_t,x,chi}
I_{t,x,\chi} := \int_{\R^{2}} e^{\iD x \cdot \xi} e^{\iD \tfrac{t}{\delta} g(\delta |\xi|)} \chi( \delta |\xi|) d\xi 
\end{equation}

We can rewrite the integral $I_{t,x,\chi}$ using  polar coordinates and the functions $h_{\pm}$ as
\begin{align*}
I_{t,x,\chi} &= \int_{\R^{+}} \int_{0}^{2\pi} e^{\iD r |x| \sin(\theta)} e^{\iD \tfrac{t}{\delta} g (\delta r)} \chi(\delta r) r d\theta dr\\
&=  \int_{\R^{+}} e^{ \iD (\tfrac{t}{\delta} g(\delta r)  + |x| r)} h_{-}( |x| r) \chi(\delta r) r dr + \int_{\R^{+}} e^{ \iD (\tfrac{t}{\delta} g(\delta r) - |x| r)} h_{+}(|x| r) \chi(\delta r) r dr\\
&:= I_{t,x,\chi,-} + I_{t,x,\chi,+}.
\end{align*}

Our first result is a low and intermediate frequency estimate assuming that $g' \neq 0$ and some other derivatives of $g$ do not vanish on a bounded interval.

\bl\label{decay2d_lowfreq_g'not0}
Let $J \subset \R^{+}$ a bounded interval. Let $\lambda> 0$ and $l \in \N$ with $l \geq 2$. Assume that $g$ satisfies (H0) and that $g$ is $\mathcal{C}^{l}$, $|g'| \geq \lambda$ and $\displaystyle \sum_{p=2}^{l} |g^{(p)}| \geq \lambda$ on $J$. Let $\chi$ be a smooth function such that $\chi(y)=0$ for any $y \notin J$. There exists $C>0$ such that for any $t \in \R^{\ast}$ and any $\delta>0$
\[
\sup_{x \in \R^2} |I_{t,x,\chi}| \leq C \frac{\delta^{\frac{2-3l}{2l}}}{|t|^{\frac{1}{2}+\frac{1}{l}}}.
\]
\el

\br\label{decay2d_lowfreq_g'cancels}
Such configuration typically occurs when there is a coupling between a high dispersive operator and a wave operator. Note that if we do not assume that $|g'| \geq \lambda$, one can only get though Van der Corput's Lemma 
\[
\sup_{x \in \R^2} |I_{t,x,\chi}| \lesssim \frac{\delta^{\frac{1-2l}{l}}}{|t|^{\frac{1}{l}}}.
\]
\er

\begin{proof}
By a change of variables, one can assume that $\delta=1$. For $p \in \{2,\cdots,l \}$ we define the sets
\[
J_{p} := \left\{y \in J, |g^{(p)}(y)| > \frac{\lambda}{2(l-1)} \right\} \bigcap_{k=2}^{p-1} \left\{ y \in J, |g^{(k)}(y)| \leq \frac{\lambda}{2(l-1)} \right\}
\]
so that by assumption $J_{p}$ is a finite union of intervals and $\sqcup_{p=2}^{l} J_{p} = J$, and the integrals
\[
I_{p,\pm} := \int_{\R^{+}} e^{ \iD (t g(r) \mp |x| r)} h_{\pm}( |x| r) \chi(r) r \mathds{1}_{J_{p}}(r) dr.
\]
We then introduce the phase $\phi_{\pm}(r)= t g(r) \mp |x| r$. We consider two cases.

\medskip
\noindent $\bullet$ Case 1: $|x| \geq \tfrac{\lambda}{2} |t|$
\medskip

Noticing that $|\phi_{\pm}^{(p)}(r)| =   |t g^{(p)}(r)| \geq \tfrac{\lambda}{l-1} |t|$, by Van der Corput's Lemma, the properties on the functions $h_{\pm}$ and since $|I_{p,\pm}| \lesssim 1$, we get
\[
|I_{p,\pm}| \lesssim \min \left( 1 , \frac{1}{|t|^{\frac{1}{p}}}  \frac{1}{\sqrt{|x|}} \right) \lesssim \min \left( 1 , \frac{1}{|t|^{\frac{1}{2}+\frac{1}{p}}} \right) \lesssim \frac{1}{|t|^{\frac{1}{2}+\frac{1}{l}}}.  
\]

\medskip
\noindent $\bullet$ Case 2: $|x| \leq \tfrac{\lambda}{2} |t|$
\medskip

Noticing that $|\phi_{\pm}'(r)| \geq \frac{\lambda}{2} |t|$, using Corollary \ref{cor_vandercorput} (note $\phi'_{\pm}$ has a finite number of zeros on $J_{p}$ since $\phi_{\pm}^{(p-1)}$ is monotonic), the properties on the functions $h_{\pm}$ and since $|I_{p,\pm}| \lesssim 1$, we get
\[
|I_{p,\pm}| \lesssim \min \left( 1 , \frac{1}{|t|} \right) \lesssim \frac{1}{(1+|t|)^{\frac{1}{2}+\frac{1}{l}}} \lesssim \frac{1}{|t|^{\frac{1}{2}+\frac{1}{l}}}.  
\]
\end{proof}

Then we provide a low frequency estimate with a better decay compare to the previous lemma when $g$ is defined and smooth close to $0$ and $g'(0) \neq 0$.

\bl\label{decay2d_ww_lowfreq}
Let $\alpha \geq 1$. Let $\Lambda,\lambda, y_{0} >0$. Assume that $g$ satisfies (H0) and is defined and $C^{1}$ in the neighborhood of $0$, that $|g'(y)| \geq \lambda$ and that $|g''(y)| \geq \lambda y^{\alpha}$ for any $y \in [0,y_{0}]$. If $\alpha=1$, assume also that $\lambda y^{\alpha+1} \leq |g'(y)-g'(0)| \leq \Lambda y^{\alpha+1}$ for any $y \in [0,y_{0}]$. Let $\chi$ be a smooth function such that $\chi(y)=0$ for any $y \geq y_{0}$. There exists $C>0$ such that for any $t \in \R^{\ast}$ and any $\delta>0$
\[
\sup_{x \in \R^2}  |I_{t,x,\chi}| \leq C \frac{\delta^{-\tfrac32\tfrac{1+\alpha}{2+\alpha}}}{|t|^{\tfrac{5+\alpha}{2(2+\alpha)}}}.
\]
\el

\br 
This lemma is typically adjusted for phases such that $g'(0) \neq 0$, $g^{(p)}(0) = 0$ and $g^{(l)}(0) \neq 0$ for $p,l \in \N$ with $2 \leq p < l$. As we will see later, the water wave phase and most of the abcd-Boussinesq phases satisfy this assumption with $l=3$ (and then $\alpha=1$).
\er

\begin{proof}
By a change of variables, one can assume that $\delta=1$.  We consider two cases.

\medskip
\noindent $\bullet$ Case 1: $|x| \leq \frac{\lambda}{2} |t|$
\medskip

Defining $\phi_{\pm}(r) =  t g(r) \mp |x| r$, we have in that case $|\phi_{\pm}'(r)| \geq \frac{\lambda}{2} |t|$ so that, using Corollary \ref{cor_vandercorput}, the properties on the functions $h_{\pm}$ and since $[0,y_{0}]$ is a bounded,
\[
|I_{t,x,\chi,\pm}| \lesssim \min \left( 1 ,\frac{1}{|t|} \right) \lesssim \frac{1}{|t|^{\tfrac{5+\alpha}{2(2+\alpha)}}}.
\]

\medskip
\noindent $\bullet$ Case 2: $|x| \geq \frac{\lambda}{2} |t|$
\medskip

We notice that $I_{t,x,\chi,\pm}$  is the evaluation of a Fourier transform (with respect to the variable $r$) at $\pm |x| - g'(0)t$
\[
I_{t,x,\chi,\pm} = \sqrt{2\pi} \mathcal{F}_{r} \left( e^{\iD t (g(r)- g'(0) r)} h_{\pm}( |x| r) \mathds{1}_{\{ r>0 \}}(r) \chi (r) r \right) (\pm |x|-g'(0)t) 
\]
so that introducing a 1d Littlewood-Paley decomposition (see Subsection \ref{ss:LittlewoodPaley}) and integrating by parts
\begin{align*}
|I_{t,x,\chi,\pm}| &\lesssim \sup_{y \in \R} \sum_{k \in \Z} \left| \int_{\R^{+}} \left( e^{-\iD r y} e^{ \iD t (g(r)- g'(0) r)} Q_{k}(r) r^{\frac12} \right) h_{\pm}( |x| r) \chi(r) r^{\frac12} dr \right|,\\
&\hspace{-0,9cm} \lesssim \sup_{y \in \R} \sum_{2^{k} \leq 2 y_{0}} \sup_{z \in [0,y_{0}]} \left| \int_{0}^{z} e^{-\iD r y} e^{ \iD t (g(r)- g'(0) r)} Q_{k}(r) r^{\frac12} dr \right| \int_{2^{k-1}}^{2^{k+1}} \left| \frac{\dD}{\dD r} \left( h_{\pm}( |x| r) \chi(r) r^{\frac12} \right) \right| dr.
\end{align*}
Using the properties on the functions $h_{\pm}$ and an easy adaptation of Lemma \ref{decay1d_VdC2} 
\[
|I_{t,x,\chi,\pm}| \lesssim \frac{1}{|t|^{\tfrac{3}{2(2+\alpha)}}} \frac{1}{\sqrt{|x|}}
\]
and the result follows in that case.
\end{proof}

Secondly, we provide estimates for more general phases.  For this reason we localize in the annulus $\mathcal{C}(2^{k-1},2^{k+1})$. Again we fix $P$ a smooth even function supported in $[-2,-\tfrac12] \cup [\frac12,2]$. We introduce the quantity for any $t \in \R$, $x \in \R^{2}$, $s \in \R$ and $k \in \mathbb{Z}$,
\[
I^{s}_{t,x,k} := \int_{\R^2} e^{\iD x \cdot \xi} e^{\iD \tfrac{t}{\delta} g(\delta |\xi|)} P(\tfrac{\delta |\xi|}{2^{k}}) |\xi|^{s} d\xi.
\]
Using a change of variables and polar coordinates we can rewrite $I_{t,x,k}$
\begin{align*}
I^{s}_{t,x,k} &= \frac{2^{(2+s)k}}{\delta^{s+2}} \int_{\frac12}^{2} \int_{0}^{2\pi} e^{\iD \frac{2^{k}}{\delta} r |x| \sin(\theta)} e^{\iD \tfrac{t}{\delta} g(2^{k} r)} P(r) r^{1+s} d\theta dr\\
&=  \frac{2^{(2+s)k}}{\delta^{s+2}} \int_{\frac12}^{2} e^{ \iD (\tfrac{t}{\delta} g(2^{k} r) + \frac{2^{k}}{\delta} |x| r)} h_{-}(\tfrac{2^{k}}{\delta} |x| r) P(r) r^{1+s} dr\\
&\hspace{1cm} + \frac{2^{(2+s)k}}{\delta^{s+2}} \int_{\frac12}^{2} e^{ \iD (\tfrac{t}{\delta} g(\delta 2^{k} r) - \frac{2^{k}}{\delta} |x| r)} h_{+}(\tfrac{2^{k}}{\delta} |x| r) P(r) r^{1+s} dr\\
&:= I^{s}_{t,x,k,-} + I^{s}_{t,x,k,+}.
\end{align*}
We recall that $\Z_{J} := \{ k \in \mathbb{Z} \text{  , } [2^{k-1},2^{k+1}] \subset J \}$.

\bl\label{decay2d_paley_1/t}
Let $J \subset \R_{+}$ an interval. Let $\lambda > 0$, $\alpha,\beta,s \in \R$. Assume that $g$ satisfies (H0), that $|g'(y)| \geq \lambda y^{\beta}$ and that $|g''(y)| \geq \lambda y^{\alpha}$ on $J$. There exists $C>0$ such that for any $t \in \R^{\ast}$, any $\delta > 0$ and any $k \in \mathbb{Z}_{J}$
\[
\sup_{x \in \R^2} |I^{s}_{t,x,k}| \leq C \frac{2^{(s+\frac{1-\alpha-\beta}{2}) k} \delta^{-(s+1)}}{|t|} + C \frac{2^{(s+1-\beta) k} \delta^{-(s+1)}}{|t|},
\]
and, if $\beta=\alpha+1$ and $l \in \R$ with $l \geq 1$, 
\[
\sup_{x \in \R^2} |I^{s}_{t,x,k}| \leq C \frac{2^{(s+\frac{2l-2-\alpha}{l}) k} \delta^{-(s+2-\frac{1}{l})}}{|t|^{\frac{1}{l}}}.
\]
\el

\begin{proof}
The second bound can be obtained by interpolation between the first bound and the fact that $|I^{s}_{t,x,k}| \lesssim \delta^{-(s+2)} 2^{(s+2)k}$. We introduce the phase $\phi_{\pm}(r)= \tfrac{t}{\delta} g(2^{k} r) \mp \tfrac{2^{k}}{\delta}  |x|r$. Let $r_{1} \in [\tfrac12,2]$ be such that  $|g'(2^{k} r_{1})| = \displaystyle \min_{r \in [1/2,2]} |g'(2^{k} r)| \gtrsim 2^{\beta k}$. We consider two different cases.

\medskip
\noindent $\bullet$ Case 1: $|x| \geq \tfrac{|t|}{2} |g'(2^{k} r_{1}) |$
\medskip

Noticing that $|\phi_{\pm}''(r)| =  2^{2k} \delta^{-1} |t g''(2^{k} r)|$, we use Van der Corput's Lemma and the properties on the functions $h_{\pm}$, we get
\[
|I^{s}_{t,x,k,\pm}| \lesssim \frac{2^{(s+2)k}}{\delta^{s+2} \sqrt{2^{2k} |\delta^{-1} t g''(2^{k} r_{2})|} }  \frac{1}{\sqrt{2^{k} \delta^{-1} |x|}} \lesssim \frac{2^{(s+2)k} \delta^{-(s+1)}}{\sqrt{2^{3k} |t g'(2^{k} r_{1}) t g''(2^{k} r_{2})|}}
\]
for some $r_{2} \in [\tfrac12,2]$. The first bound follows in that case.

\medskip
\noindent $\bullet$ Case 2: $|x| \leq \tfrac{|t|}{2} \displaystyle |g'(2^{k} r_{1})|$
\medskip

Noticing that $|\phi_{\pm}'(r)| \gtrsim 2^{k} \delta^{-1} |t g'(2^{k} r_{1})|$ and using Corollary \ref{cor_vandercorput} and the properties on the functions $h_{\pm}$, we get
\[
|I^{s}_{t,x,k,\pm}| \lesssim \frac{2^{(s+2)k}}{\delta^{s+2}} \frac{1}{2^{k} \delta^{-1} |t g'(2^{k} r_{1})|}.
\]
The first bound follows in that case.
\end{proof}

The next two lemmas can be seen as a generalization of the Van der Corput's Lemma to the 2d case. 

\bl\label{VdC2d}
Let $J \subset \R^{\ast}_{+}$ an interval. Let $\Lambda,\lambda> 0$, $\alpha,s \in \R$ with $\alpha \neq -2$. Assume that $g$ satisfies (H0), that $|g'(y)| \geq \lambda y^{\alpha+1}$ and that $|g''(y)| \geq \lambda y^{\alpha}$ on $J$. 
\begin{enumerate}
\item If $\alpha > s > -2$ or $\alpha <s < -2$, there exists $C>0$ such that for any $t \in \R^{\ast}$ and any $\delta>0$
\[
\sum_{k \in \Z_{J}} \sup_{x \in \R^{2}} \left| I^{s}_{t,x,k} \right| \leq C \frac{\delta^{-\frac{(\alpha+1)(s+2)}{2+\alpha}}}{|t|^{\frac{s+2}{2+\alpha}}} \text{     ,     }  \sup_{x \in \R^{2}} \left| \int_{\R^2} e^{\iD x \cdot \xi} e^{\iD \tfrac{t}{\delta} g(\delta |\xi|)} \sum_{k \in \Z_{J}} Q_{k}(\delta |\xi|) |\xi|^{s} d\xi \right| \leq C \frac{\delta^{-\frac{(\alpha+1)(s+2)}{2+\alpha}}}{|t|^{\frac{s+2}{2+\alpha}}}.
\]
\item If $s=\alpha$ and $\alpha \notin \{ -2 , -1 \}$, assuming furthermore that $g$ is $\mathcal{C}^{3}$, that $|g'(y)| \leq \Lambda y^{\alpha+1}$, that $|g''(y)| \leq \Lambda y^{\alpha}$ and that $|g'''(y)| \leq \Lambda y^{\alpha-1}$ on $J$, there exists $C>0$ such that for any $t \in \R^{\ast}$ and any $\delta>0$
\[
\sup_{x \in \R^{2}} \;\; \sum_{k \in \Z_{J}}\left| I^{\alpha}_{t,x,k} \right| \leq C  \frac{\delta^{-(\alpha+1)}}{|t|} \text{     ,     } \sup_{x \in \R^{2}} \left| \int_{\R^2} e^{\iD x \cdot \xi} e^{\iD \tfrac{t}{\delta} g(\delta |\xi|)} \sum_{k \in \Z_{J}} Q_{k}(\delta |\xi|) |\xi|^{\alpha} d\xi \right| \leq C  \frac{\delta^{-(\alpha+1)}}{|t|}.
\]
\end{enumerate}
\el

\br\label{exp_VdC2d}
When $J=]0,1]$, $\delta=1$ and $s=0$, we get item (c) of Theorem 1 in \cite{decay_estimates_nonhomogeneous_guo_et_al}. 

This lemma can be read as follows. We assume by simplicity that $\delta=1$ and $s=0$. Let $\Lambda,\lambda>0$, $l \in \mathbb{N}$ with $l \geq 2$. We recall that we define the quantity $I_{t,x,\chi}$ in \eqref{I_t,x,chi}.
\begin{enumerate}
\item Assume that $|g^{(p)}| \geq \lambda$ on a bounded interval $J \subset \R^{+}$ containing $0$ with $g^{(k)}(0)=0$ for any $k \in \{1,\cdots,l-1 \}$, then for any smooth function $\chi$ whose support is a subset of $J$
\[
| I_{t,x,\chi} | \lesssim \frac{1}{|t|^{\frac{2}{p}}}.
\]
\item Assume that $\Lambda \geq |g^{(p)}| \geq \lambda$ on an unbounded interval $J \subset \R^{+}$ with $\inf(J)>0$, then there exists $y_{1} \in J$ such that for any smooth function $\chi$ whose support is a subset of $[0,y_{1}]$
\[
| I_{t,x,1-\chi} | \lesssim \frac{1}{|t|^{\frac{2}{p}}}.
\]
\end{enumerate}
Note that when $J=(0,y_{0}]$ and $\alpha <s < -2$ or $J=[y_{0},\infty)$ and $\alpha > s > -2$ for some $y_{0}>0$, we rather use Lemma \ref{decay2d_paley_1/t} since the decay provided by Lemma \ref{VdC2d} is not optimal.

Finally, as noted in Remark \ref{r:alphanot2} about Lemma \ref{decay1d_VdC2}, the assumption $\alpha \neq -2$ is crucial, the fact that $(\alpha,s) \neq (-1,-1)$ is due to our proof and we interpret this issue as the fact that $\alpha=-1$ corresponds to a wave-type behavior so that the dispersive effects can not be optimal.
\er 

\begin{proof}
We begin with the first estimate. In that case $(s+2)(s-\alpha)<0$. We have, using Lemma \ref{decay2d_paley_1/t}, for any $k \in \Z_{J}$
\[
|I^{s}_{t,x,k}| \lesssim \min \left( \frac{2^{(2+s)k}}{\delta^{s+2}} , \frac{2^{(s-\alpha) k}}{\delta^{s+1} |t|} \right) 
\]
so that if $\alpha>-2$
\[
\sum_{k \in \Z_{J}} \sup_{x \in \R^{2}} \left| I^{s}_{t,x,k} \right| \lesssim \sum_{2^{k} |t|^{\frac{1}{\alpha+2}} \leq \delta^{\frac{1}{\alpha+2}}} \frac{2^{(2+s)k}}{\delta^{s+2}} + \sum_{2^{k} |t|^{\frac{1}{\alpha+2}} \geq \delta^{\frac{1}{\alpha+2}}} \frac{2^{(s-\alpha) k}}{\delta^{s+1} |t|} ,
\]
whereas if $\alpha<-2$
\[
\sum_{k \in \Z_{J}} \sup_{x \in \R^{2}} \left| I^{s}_{t,x,k} \right| \lesssim \sum_{2^{k} |t|^{\frac{1}{\alpha+2}} \leq \delta^{\frac{1}{\alpha+2}}} \frac{2^{(s-\alpha) k}}{\delta^{s+1} |t|} + \sum_{2^{k} |t|^{\frac{1}{\alpha+2}} \geq \delta^{\frac{1}{\alpha+2}}} \frac{2^{(2+s)k}}{\delta^{s+2}}.
\]
The first estimate follows and also the second one by taking $P=Q_{0}$ in the definition of $I^{s}_{t,x,k}$. Then we consider the third estimate and define the phase $\phi_{\pm}(r)= \tfrac{t}{\delta} g(2^{k}r) \mp \frac{2^{k}}{\delta} |x| r$. Two cases occurs.

\medskip
\noindent $\bullet$ Case 1: $\frac{\lambda}{4} 2^{(\alpha+1) k} |t| \leq |x| \leq 4 \Lambda 2^{(\alpha+1) k} |t|$
\medskip

Here $2^{k} |x| \sim 2^{(\alpha+2)k} |t|$ and using Van der Corput's Lemma and the properties on $h_{\pm}$
\[
|I^{\alpha}_{t,x,k}| \lesssim \frac{2^{(\alpha+2)k}}{\delta^{\alpha+2} \sqrt{2^{(\alpha+2)k} \delta^{-1} |t|}} \frac{1}{\sqrt{\frac{2^{k}}{\delta} |x|}} \lesssim \frac{\delta^{-(\alpha+1)}}{|t|}.
\]

Note that the set $A_{t,x} := \{ k \in \mathbb{Z}_{J}, \frac{\lambda}{4} 2^{(\alpha+1) k} |t| \leq |x| \leq 4 \Lambda 2^{(\alpha+1) k} |t| \}$ is bounded by a number independent of $t,x,\delta$.

\medskip
\noindent $\bullet$ Case 2:  $\frac{\lambda}{4} 2^{(\alpha+1) k} |t| \geq |x|$ or $|x| \geq 4 \Lambda 2^{(\alpha+1) k} |t|$
\medskip

We notice that $|\phi'_{\pm}(r)| \geq \frac{\lambda}{4} \delta^{-1} 2^{(\alpha+2) k} |t|$. Integrating by parts we get
\[
I^{\alpha}_{t,x,k,\pm} = \iD \frac{2^{(2+\alpha)k}}{\delta^{\alpha+2}} \int_{\frac12}^{2} e^{ \iD \phi_{\pm}'(r)} \frac{d}{dr} \left( \frac{1}{\phi_{\pm}'(r)} h_{\pm}(\tfrac{2^{k}}{\delta} |x| r) P(r) r \right) dr
\]
so that by Corollary \ref{cor_vandercorput}, the properties on $h_{\pm}$ and the bounds on $g',g'',g'''$
\[
|I^{\alpha}_{t,x,k,\pm}| \lesssim \frac{2^{-(\alpha+2) k} \delta^{-\alpha}}{t^2}.
\]

Therefore gathering all previous estimates, we get if $\alpha>-2$ and $\alpha \neq -1$
\begin{align*}
\sum_{k \in \Z_{J}} \left|  I^{\alpha}_{t,x,k} \right| &\lesssim \sum_{k \in A_{t,x}} | I^{\frac{\alpha}{2}}_{t,x,k}| + \sum_{k \in \Z_{J} \text{, } k \notin A_{t,k}} | I^{\frac{\alpha}{2}}_{t,x,k}|\\
&\lesssim \frac{\delta^{-(\alpha+1)}}{|t|} + \sum_{2^{k} |t|^{\frac{1}{\alpha+2}} \leq \delta^{\frac{1}{\alpha+2}}} \frac{2^{(\alpha+2) k}}{\delta^{\alpha+2}} + \sum_{2^{k} |t|^{\frac{1}{\alpha+2}} \geq |t|^{\frac{1}{\alpha+2}}} \frac{2^{-(\alpha+2) k} \delta^{-\alpha}}{t^2},
\end{align*}
whereas if $\alpha <-2$
\[
\sum_{k \in \Z_{J}} \left|  I^{\alpha}_{t,x,k} \right| \lesssim \frac{\delta^{-(\alpha+1)}}{|t|} + \sum_{2^{k} |t|^{\frac{1}{\alpha+2}} \leq |t|^{\frac{1}{\alpha+2}}} \frac{2^{-(\alpha+2) k} \delta^{-\alpha}}{t^2} + \sum_{2^{k} |t|^{\frac{1}{\alpha+2}} \geq |t|^{\frac{1}{\alpha+2}}} \frac{2^{(\alpha+2) k}}{\delta^{\alpha+2}}.
\]
The third and fourth estimates follow.
\end{proof}

When $|g'(y)| \gtrsim 1$ and $|g''(y)| \gtrsim y^{\alpha}$ with $\alpha<-1$ and for any $y \geq y_{1} > 0$, Lemma \ref{decay2d_paley_1/t} provides a decay of order $\tfrac{1}{|t|}$ whereas Lemma \ref{VdC2d} does not. We studied a similar situation in the low frequency case in Lemma \ref{decay2d_ww_lowfreq}. The next lemma provides bounds in the high frequency case.

\bl\label{decay2d_theta=0_HF}
Let $\Lambda,\lambda,y_{1}>0$, $\ell \in \R^{\ast}$ and $\alpha \in \R$ with $\alpha<-1$. Assume that $g$ satisfies (H0), that $|g'(y)| \geq \lambda$ and that $|g''(y)| \geq \lambda y^{\alpha}$ for any $y \geq y_{1}$. 
\begin{enumerate}
\item If $\frac{\alpha-1}{2} < s < -2$ or $\frac{\alpha-1}{2} > s > -2$, there exists $C>0$ such that then for any $t \in \R^{\ast}$ and any $\delta>0$
\[
\sum_{2^{k} \geq 2 y_{1}} \sup_{x \in \R^{2}} \left| I^{s}_{t,x,k} \right| \leq C  \frac{\delta^{-\frac{(\alpha+1)(s+2)}{\alpha+3}}}{|t|^{\frac{2(s+2)}{\alpha+3}}} \text{     ,     } \sup_{x \in \R^2} \left| \int_{\R^2} e^{\iD x \cdot \xi} e^{\iD \tfrac{t}{\delta} g(\delta |\xi|)} \sum_{2^{k} \geq 2 y_{1}} Q_{k}(\delta \xi)  |\xi|^{s} d\xi \right|  \leq C  \frac{\delta^{-\frac{(\alpha+1)(s+2)}{\alpha+3}}}{|t|^{\frac{2(s+2)}{\alpha+3}}}.
\]
\item If furthermore $\alpha \neq -2$ and $\lambda y^{\alpha+1} \leq |g'(y)-\ell| \leq \Lambda y^{\alpha+1}$ for any $y \geq y_{1}$, there exists $C>0$ such that for any $t \in \R^{\ast}$ and any $\delta>0$
\[
\sup_{x \in \R^{2}} \;\; \sum_{2^{k} \geq 2 y_{1}} \left| I^{\frac{\alpha-1}{2}}_{t,x,k} \right| \leq C  \frac{\delta^{-\frac{\alpha+1}{2}}}{|t|} \text{     ,     } \sup_{x \in \R^2} \left| \int_{\R^2} e^{\iD x \cdot \xi} e^{\iD \tfrac{t}{\delta} g(\delta |\xi|)} \sum_{2^{k} \geq 2 y_{1}} Q_{k}(\delta \xi)  |\xi|^{\frac{\alpha-1}{2}} d\xi \right| \leq C  \frac{\delta^{-\frac{\alpha+1}{2}}}{|t|}.
\]
\end{enumerate}
\el

\br 
This kind of configuration occurs for some abcd-Boussinesq systems, see Subsection \ref{ss:abcd_Boussi}.
\er

\begin{proof}

We begin with the first estimate. In that case $(s+2)(s-\frac{\alpha-1}{2})<0$. Using Lemma \ref{decay2d_paley_1/t}, for any $2^{k} \delta \geq 2 y_{1}$
\[
|I^{s}_{t,x,k}| \lesssim \min \left( \frac{2^{(2+s)k}}{\delta^{s+2}} , \frac{2^{(s+\frac{1-\alpha}{2}) k} \delta^{-(s+1)}}{|t|} \right) 
\]
so that if $\alpha>-3$
\[
\sum_{2^{k} \geq 2 y_{1}} \sup_{x \in \R^{2}} \left| I^{s}_{t,x,k} \right| \lesssim \sum_{2^{k} |t|^{\frac{2}{\alpha+3}} \leq \delta^{\frac{2}{\alpha+3}}} \frac{2^{(2+s)k}}{\delta^{s+2}} + \sum_{2^{k} |t|^{\frac{2}{\alpha+3}} \geq \delta^{\frac{2}{\alpha+3}}} \frac{2^{(s+\frac{1-\alpha}{2}) k} \delta^{-(s+1)}}{|t|},
\]
whereas if $\alpha<-3$
\[
\sum_{2^{k} \geq 2 y_{1}} \sup_{x \in \R^{2}} \left| I^{s}_{t,x,k} \right| \lesssim \sum_{2^{k} |t|^{\frac{2}{\alpha+3}} \leq \delta^{\frac{2}{\alpha+3}}} \frac{2^{(s+\frac{1-\alpha}{2}) k} \delta^{-(s+1)}}{|t|} + \sum_{2^{k} |t|^{\frac{2}{\alpha+3}} \geq \delta^{\frac{2}{\alpha+3}}} \frac{2^{(2+s)k}}{\delta^{s+2}}.
\]
The first estimate follows and also the second one by taking $P=Q_{0}$ in the definition of $I^{s}_{t,x,k}$. Then we consider the third estimate. Two cases occurs.

\medskip
\noindent $\bullet$ Case 1: $|x| \leq \frac{\lambda}{2} |t|$
\medskip

Defining $\phi_{\pm}(r) =  \tfrac{t}{\delta} g(2^{k} r) \mp \frac{2^{k}}{\delta} |x| r$, we have in that case $|\phi_{\pm}'(r)| \geq \frac{\lambda}{2} 2^{k} \delta^{-1} |t|$ so that, using Corollary \ref{cor_vandercorput} and the properties on the functions $h_{\pm}$
\[
|I^{\frac{\alpha-1}{2}}_{t,x,k}| \lesssim \frac{2^{(\frac{\alpha+1}{2})k} \delta^{-\frac{\alpha+1}{2}}}{|t|}
\]
and since $\alpha+1<0$, the result follows in that case.

\medskip
\noindent $\bullet$ Case 2: $|x| \geq \frac{\lambda}{2} |t|$
\medskip

We notice that $I^{\frac{\alpha-1}{2}}_{t,x,k,\pm}$  is the evaluation of a Fourier transform (with respect to the variable $r$) at $\pm |x| - \ell t$
\[
I^{\frac{\alpha-1}{2}}_{t,x,k,\pm}= \sqrt{2 \pi} \mathcal{F}_{r} \left( e^{\iD \tfrac{t}{\delta} (g(\delta r)- \ell \delta r)} P(\tfrac{\delta}{2^{k}} r) \mathds{1}_{\{ r>0 \}}(\delta r) h_{\pm}( |x| r) r^{\frac{\alpha+1}{2}} \right) (\pm |x|-\ell t).
\]
Denoting $y_{\pm}:=\pm |x|-\ell t$, integrating by parts and using the properties on the functions $h_{\pm}$ we obtain
\begin{align*}
|I^{\frac{\alpha-1}{2}}_{t,x,k,\pm}| &\lesssim \left| \int_{\R^{+}} \left( e^{-\iD r y_{\pm}} e^{ \iD \tfrac{t}{\delta} (g(\delta r)- \ell \delta r)} P(\tfrac{\delta}{2^{k}} r) r^{\frac{\alpha}{2}} \right) h_{\pm}( |x| r) r^{\frac12} dr \right|,\\
&\hspace{-1cm}\lesssim \sup_{z \in [2^{k-1},2^{k+1}]} \left| \int_{2^{k-1}}^{z} e^{-\iD r y_{\pm}} e^{ \iD \tfrac{t}{\delta} (g(\delta r)- \ell \delta r)} P(\tfrac{\delta}{2^{k}} r) r^{\frac{\alpha}{2}} dr \right| \left(  \frac{1}{\sqrt{|x|}} +\int_{2^{k-1}}^{2^{k+1}} \left| \frac{\dD}{\dD r} \left( h_{\pm}( |x| r) r^{\frac12} \right) \right| \right) dr,\\
&\hspace{-1cm}\lesssim \sup_{z \in [2^{k-1},2^{k+1}]} \left| \int_{2^{k-1}}^{z} e^{-\iD r y_{\pm}} e^{ \iD \tfrac{t}{\delta} (g(\delta r)- \ell \delta r)} P(\tfrac{\delta}{2^{k}} r) r^{\frac{\alpha}{2}} dr \right| \frac{1}{\sqrt{|t|}}.
\end{align*}
The result in that case follows from an easy adaptation of Lemma \ref{decay1d_VdC2}.
\end{proof}

We then give two other lemmas with weaker decays compare to the previous lemma. They provide however a decay when $g'$ has a positive zero and the bounds are uniform with respect $\delta \in (0,1]$. We see the problem as a perturbation of the half-wave equation when $\delta$ is small. The following lemma gives a decay when $g'$ and $g''$ do not vanish at the same time.

\bl\label{decay2d_wave_2deriv}
 Let $\alpha \in \R$, $\lambda >0$ and $y_{1} \geq y_{0} > 0$. Assume that $g$ satisfies (H0). 
\begin{enumerate}
\item If $|g'| \geq \lambda$ and $g''$ has a finite number of zeros on $(0,y_{0}]$, there exists $C>0$ such that for any $t \in \R^{\ast}$, any $\delta>0$ and any $k \in \mathbb{Z}$ such that $2^{k} \leq \frac12 y_{0}$ 
\[
\sup_{x \in \R^2} |I^{s}_{t,x,k}| \leq C \frac{2^{(s+\frac32) k} \delta^{-(s+\frac32)}}{\sqrt{|t|}}.
\]
\item If $|g'| + |g''| \geq \lambda$ and $g''$ has a finite number of zeros on $[\tfrac{1}{2} y_{0} , 2 y_{1} ]$, there exists $C>0$ such that for any $t \in \R^{\ast}$, any $\delta>0$ and any $k \in \mathbb{Z}$ such that $y_{0} \leq 2^{k} \leq y_{1}$
\[
\sup_{x \in \R^2} |I^{s}_{t,x,k}| \leq C \frac{\delta^{-(s+\frac32)}}{\sqrt{|t|}} \text{     ,    } \sup_{x \in \R^2} \left| \int_{\R^2} e^{\iD x \cdot \xi} e^{\iD \tfrac{t}{\delta} g(\delta |\xi|)} \sum_{y_{0} \leq 2^{k} \leq y_{1}} Q_{k}(\delta \xi)  |\xi|^{s} d\xi \right| \leq C \frac{\delta^{-(s+\frac32)}}{\sqrt{|t|}}.
\]
\item If $|g'(y)| \geq \lambda y^{\alpha+1}$ and $|g''(y)| \geq \lambda y^{\alpha}$ for any $y \geq y_{1}$ with $\alpha \neq -2$, there exists $C>0$ such that for any $t \in \R^{\ast}$ and any $\delta>0$
\[
\sum_{2^{k} \geq 2 y_{1}} \sup_{x \in \R^{2}} |I^{\frac{\alpha}{2}-1}_{t,x,k}| \leq C \frac{\delta^{-\frac{\alpha+1}{2}}}{\sqrt{|t|}} \text{     ,     } \sup_{x \in \R^2} \left| \int_{\R^2} e^{\iD x \cdot \xi} e^{\iD \tfrac{t}{\delta} g(\delta |\xi|)} \sum_{2^{k} \geq 2 y_{1}} Q_{k}(\delta \xi)  |\xi|^{\frac{\alpha}{2}-1} d\xi \right| \leq C \frac{\delta^{-\frac{\alpha+1}{2}}}{\sqrt{|t|}}.
\]
\item  If $|g'(y)| \geq \lambda $ and $|g''(y)| \geq \lambda y^{\alpha}$ for any $y \geq y_{1}$ with $\alpha < -1$ and $\alpha \neq -3$, there exists $C>0$ such that for any $t \in \R^{\ast}$ and any $\delta>0$
\[
\sum_{2^{k} \geq 2 y_{1}} \sup_{x \in \R^{2}} |I^{\frac{\alpha-5}{4} }_{t,x,k}| \leq C \frac{\delta^{-\frac{\alpha+1}{2}}}{\sqrt{|t|}} \text{     ,     } \sup_{x \in \R^2} \left| \int_{\R^2} e^{\iD x \cdot \xi} e^{\iD \tfrac{t}{\delta} g(\delta |\xi|)} \sum_{2^{k} \geq 2 y_{1}} Q_{k}(\delta \xi)  |\xi|^{\frac{\alpha-5}{4}} d\xi \right| \leq C \frac{\delta^{-\frac{\alpha+1}{2}}}{\sqrt{|t|}}.
\]
\end{enumerate} 
\el
\br 
By taking $s=-\frac{3}{2}$ in the first two estimates these bounds are uniform with respect to $\delta \in (0,1]$ and the case $\delta=0$ exactly corresponds to the bounds of the 2d wave equation, see for instance Lemma 2.1 in the Oberwolfach seminar by M. Vi\c{s}an \cite{KT_Oberwolfach}. By taking $\alpha \leq -1$ (which is always possible) in the last two estimates we also have uniform bounds with respect to $\delta \in (0,1]$.
\er

\begin{proof}
We begin with items (1) and (2). We introduce the following quantities to bound
\begin{align*}
&I_{1} : = \frac{2^{(s+2)k}}{\delta^{s+2}} \int_{\frac12}^{2} e^{ \iD (\sgn(t g') \frac{2^{k}}{\delta} |x| r + \tfrac{t}{\delta} g(2^{k} r))} \mathds{1}_{\{ |g'| > \tfrac{\lambda}{4} \}}(2^{k} r)  P(r) r^{s+1} J_{0}( \tfrac{2^{k}}{\delta} |x| r) e^{-\iD \; \sgn(t g') \frac{2^{k}}{\delta} |x| r} dr,\\
&I_{2} :=  \frac{2^{(s+2)k}}{\delta^{s+2}} \int_{0}^{2\pi} \int_{\frac12}^{2} e^{ \iD (\frac{2^{k}}{\delta} |x| r \sin(\theta) + \tfrac{t}{\delta} g(2^{k} r))} \mathds{1}_{\{ |g''| > \tfrac{\lambda}{4} , |g'| \leq \tfrac{\lambda}{4} \}}(2^{k} r) P(r) r^{s+1} dr d\theta.
\end{align*}
Note that the sets $\{ |g'| > \tfrac{\lambda}{4} \}$ and $\{ |g''| > \tfrac{\lambda}{4} , |g'| \leq \tfrac{\lambda}{4} \}$ are finite unions of intervals. Defining $\phi_{1}(r)=\sgn(t g') \tfrac{2^{k}}{\delta} |x| r + \tfrac{t}{\delta} g(2^{k} r)$ so that $|\phi_{1}'(r)| =2^{k} \delta^{-1} |x| + 2^{k} \delta^{-1} |t g'(2^{k} r)| \gtrsim 2^{k} \delta^{-1} (|x| + |t|)$
and using Corollary \ref{cor_vandercorput} and some properties of the Bessel function $J_{0}$ we obtain 
\begin{align*}
|I_{1}| &\lesssim \min \left( \frac{2^{(s+2)k}}{\delta^{s+2}} , \frac{2^{(s+2)k} \delta^{-(s+2)}}{ \frac{2^{k}}{\delta} |x| + \frac{2^{k}}{\delta} |t| }  \sqrt{1+ \tfrac{2^{k}}{\delta} |x|} \right)\\
&\lesssim \min \left( \frac{2^{(s+2)k}}{\delta^{s+2}} , \frac{2^{(s+1)k} \delta^{-(s+1)}}{|t|} \mathds{1}_{\{ \frac{2^{k}}{\delta} |t| \leq 2 \}} + \frac{2^{(s+\frac32) k} \delta^{-(s+\frac32)}}{\sqrt{|t|}} \mathds{1}_{\{ \frac{2^{k}}{\delta} |t| \geq 2 \}} \right).
\end{align*}
It provides the desired bound for $I_{1}$. Then, we introduce $\phi_{2}(r)=\frac{2^{k}}{\delta} |x| r \sin(\theta) + \tfrac{t}{\delta} g(2^{k} r)$ so that $|\phi_{2}''(r)| \gtrsim 2^{2k} \delta^{-1} |t|$ and Van der Corput's Lemma provides the desired bound for $I_{2}$ since we only estimate $I_{2}$ when $2^{k} \sim 1$. The third bound follows by boundedness of the domain of summation.
Item (3) follows from Lemma \ref{VdC2d} with $s=\frac{\alpha}{2}-1$ whereas item (4) is a consequence of Lemma \ref{decay2d_theta=0_HF} with $s=\frac{\alpha-5}{4}$.
\end{proof}

Then, we consider the situation where some derivatives of $g$ do not vanish at the same time. By simplicity, we only consider low and intermediate frequencies.

\bl\label{decay2d_wave_multideriv}
Let $\lambda > 0$, $\alpha \in \R$, $y_{1} \geq y_{0} > 0$ and $l \in \N$ with $l \geq 2$. Assume that $g$ satisfies (H0).
\begin{enumerate}
\item If $|g'| \geq \lambda$ and $g''$ has a finite number of zeros on $(0,y_{0}]$, there exists $C>0$ such that for any $t \in \R^{\ast}$, any $\delta>0$ and for any $k \in \mathbb{Z}$ such that $2^{k} \delta \leq \frac12 y_{0}$
\[
\sup_{x \in \R^2} |I^{s}_{t,x,k}| \leq C \frac{2^{(s+2-\frac{1}{l})k} \delta^{-(s+2-\frac{1}{l})}}{|t|^{\frac{1}{l}}}.
\]
\item If $g$ is $\mathcal{C}^{l}$, $\displaystyle \sum_{p=1}^{l} |g^{(p)}| \geq \lambda$ and $g''$ has a finite number of zeros on $[\tfrac{1}{2} y_{0} , 2 y_{1} ]$, there exists $C>0$ such that for any $t \in \R^{\ast}$, any $\delta>0$ and for any $k \in \mathbb{Z}$ such that $y_{0} \leq 2^{k} \leq y_{1}$
\[
\sup_{x \in \R^2} |I^{s}_{t,x,k}| \leq C \frac{\delta^{-(s+2-\frac{1}{l})}}{|t|^{\frac{1}{l}}}  \text{      ,      } \sup_{x \in \R^2} \left| \int_{\R^2} e^{\iD x \cdot \xi} e^{\iD \tfrac{t}{\delta} g(\delta |\xi|)} \sum_{y_{0} \leq 2^{k} \leq y_{1}} Q_{k}(\delta \xi)  |\xi|^{s} d\xi \right| \leq C \frac{\delta^{-(s+2-\frac{1}{l})}}{|t|^{\frac{1}{l}}}.
\]
\end{enumerate} 
\el

\begin{proof}
For $p \in \{1,\dots,l \}$, we define the sets
\[
J_{p} := \left\{y \in [\tfrac12 y_{0}, 2 y_{1}], |g^{(p)}(y)| > \frac{\lambda}{2(l-1)} \right\} \bigcap_{k=1}^{p-1} \left\{ y \in [\tfrac12 y_{0}, 2 y_{1}], |g^{(k)}(y)| \leq \frac{\lambda}{2(l-1)} \right\}
\]
so that by assumption $J_{p}$ is a finite union of intervals and $\sqcup_{p=1}^{l} J_{p} = [\tfrac12 y_{0}, 2 y_{1}]$, and we define the integrals
\[
I_{p} = \frac{2^{(s+2)k}}{\delta^{(s+2)}} \int_{0}^{2\pi} \int_{\frac12}^{2} e^{ \iD ( \frac{2^{k}}{\delta} |x| r \sin(\theta) + \tfrac{t}{\delta} g(2^{k} r))} \mathds{1}_{ J_{p}}(2^{k}r) P(r) r^{s+1} dr d\theta.
\]
As in the proof of Lemma \ref{decay2d_wave_2deriv}, we have
\[
|I_{1}| \lesssim \frac{2^{(s+\frac32) k} \delta^{-(s+\frac32)}}{\sqrt{|t|}}.
\]
For $p \geq 2$, defining $\phi(r) = \frac{2^{k}}{\delta} |x| r \sin(\theta) + \tfrac{t}{\delta} g(2^{k} r)$, Van der Corput's Lemma gives
\[
|I_{p}| \lesssim\frac{2^{(s+2)k}}{\delta^{s+2}} \left( \frac{\delta}{2^{pk} |t|} \right)^{\frac{1}{p}} \lesssim 2^{k(s+1)} \frac{\delta^{-(s+\frac32)}}{\sqrt{|t|}} + 2^{(s+1)k} \frac{\delta^{-(s+2-\frac{1}{l})}}{|t|^{\frac{1}{l}}},
\]
Then for any $p \in \{1,\dots,l \}$, if $2^{k} |t| \leq \delta$ we have
\[
|I_{t,x,k}| \lesssim \frac{2^{(s+2)k}}{\delta^{s+2}} = \frac{2^{(s+2-\frac{1}{l})k}}{\delta^{s+2}} 2^{\frac{k}{l}} \lesssim 2^{(s+2-\frac{1}{l})k} \frac{\delta^{-(s+2-\frac{1}{l})}}{|t|^{\frac{1}{l}}}
\]
whereas if $2^{k} |t| \geq \delta$ we notice that
\[
\frac{\delta^{-(s+\frac32)}}{\sqrt{|t|}} = \frac{\delta^{-(s+\frac32)}}{|t|^{\frac{1}{l}}} \frac{1}{|t|^{(\frac{1}{2}-\frac{1}{l})}} \leq \frac{\delta^{-(s+\frac32)}}{|t|^{\frac{1}{l}}} (2^{k} \delta^{-1})^{\frac{1}{2}-\frac{1}{l}} = 2^{(\frac12 -\frac{1}{l})k} \frac{\delta^{-(s+2-\frac{1}{l})}}{|t|^{\frac{1}{l}}}.
\]
The first and second bounds follow using that $2^{k} \sim 1$ for the second bound. The third bound is a consequence of the boundedness of the domain of summation and by taking $P=Q_{0}$ in the definition of $I_{t,x,k}$.
\end{proof}

\subsection{Strichartz estimates}\label{ss:Strichartz}

In this subsection, we prove Strichartz estimates from the $L^{1} \to L^{\infty}$ type bounds obtained in the previous subsections.  We first recall a few facts. For $n=1$ and $n=2$ and a smooth function $P$ supported in the annulus $\mathcal{C}(\frac12,2)$ we defined
\[
I^{s}_{t,x,k} :=   \int_{\R^{+}} e^{\iD x \xi} e^{\iD \tfrac{t}{\delta} g(\delta |\xi|)} P(\tfrac{\delta \xi}{2^{k}}) \xi^{s} d\xi  \text{  or   } \int_{\R^2} e^{\iD x \cdot \xi} e^{\iD \tfrac{t}{\delta} g(\delta |\xi|)} P( \tfrac{\delta |\xi|}{2^{k}}) |\xi|^{s} d\xi.
\]
We defined a Littlewood-Paley decomposition $(Q_{j})_{j \in \Z}$ in Subsection \ref{ss:LittlewoodPaley} and we note that if $P$ is equal to $1$ in the annulus $\mathcal{C}(\tfrac35,\tfrac85)$, then for any $k \in \Z$, $Q_{k} P (\tfrac{ \cdot}{2^{k}}) = Q_{k}$.

For any $\sigma \in (0,1]$ we say that an ordered pair $(q,r)$ is sharp $\sigma$-admissible (definition introduced by \cite{KeelTao98}) if $q,r \geq 2$ and
\be\label{sharpadmissible}
\frac{1}{q} + \frac{\sigma}{r} = \frac{\sigma}{2} \text{   ,   } (q,r,\sigma) \neq (2,\infty,1).
\ee

We first state Strichartz estimates similar to the wave equation.

\bpr\label{Strichartz_Paley}
Let $n=1$ or $2$, $J \subset \R^{+}$. Assume that $g$ satisfies (H0). Suppose there exists $C>0$, $\sigma \in (0,1]$ and $\beta, \gamma \in \R$ such that for any $t \in \R^{\ast}$, any $\delta>0$ and any $k \in \Z_{J}$
\[
\sup_{x \in \R^{n}} | I_{t,x,k}^{0} | \leq C 2^{\gamma k} \frac{\delta^{\beta}}{|t|^{\sigma}}
\]
for some $P$ supported in the annulus $\mathcal{C}(\tfrac12,2)$ that is equal to $1$ in the annulus $\mathcal{C}(\tfrac35,\tfrac85)$. Then, for any $(q,\tilde q,r, \tilde r) \in [2,\infty]^4$ with $(q,r)$ and $(\tilde q,\tilde r)$ sharp $\sigma$-admissible, there exists a constant $C_{1}>0$ such that for any $\delta>0$, any $k \in \Z_{J}$, any $f \in L^{2}(\R^{n})$ and any $F \in L^{q'}_{t}(\R ; L^{r'}_{x}(\R^{n}))$
\begin{align*}
&\left\| e^{\iD \frac{t}{\delta} g(\delta |D|)}  Q_{k}(\delta |D|) f \right\|_{L^{q}_{t} L^{r}_{x}} \leq C_{1} 2^{\gamma (\frac12 - \frac{1}{r}) k} \delta^{\beta (\frac12-\frac{1}{r})}  \| Q_{k}(\delta |D|) f \|_{L^{2}},\\
&\left\| \int_{\R} e^{-\iD \frac{s}{\delta} g(\delta |D|)} Q_{k}(\delta |D|) F(s) ds \right\|_{L^{2}_{x}} \leq C_{1} 2^{\gamma (\frac12 - \frac{1}{r}) k} \delta^{\beta (\frac12-\frac{1}{r})}  \| Q_{k}(\delta |D|) F \|_{L^{q'}_{t} L^{r'}_{x}},\\
&\left\| \int_{\R}  e^{\iD \frac{(t-s)}{\delta} g(\delta |D|)} Q_{k}(\delta |D|)  F(s) ds \right\|_{L^{q}_{t} L^{r}_{x}}  \leq C_{1} 2^{\gamma (1 - \frac{1}{r} - \frac{1}{\tilde r}) k}  \delta^{\beta (1-\frac{1}{r}-\frac{1}{\tilde r})} \| Q_{k}(\delta |D|) F \|_{L^{\tilde q'}_{t} L^{\tilde r'}_{x}},\\
&\left\| \int_{s < t} e^{\iD \frac{(t-s)}{\delta} g(\delta |D|)} Q_{k}(\delta |D|)  F(s) ds \right\|_{L^{q}_{t} L^{r}_{x}}  \leq C_{1} 2^{\gamma (1 - \frac{1}{r} - \frac{1}{\tilde r}) k} \delta^{\beta (1-\frac{1}{r}-\frac{1}{\tilde r})} \| Q_{k}(\delta |D|) F \|_{L^{\tilde q'}_{t} L^{\tilde r'}_{x}}.
\end{align*}
Furthermore, if $\chi$ is a smooth bounded function such that $\chi \sum_{k \in \Z_{J}} Q_{k} = \chi$, then for any $(q,\tilde q,r, \tilde r) \in [2,\infty]^4$ with $(q,r)$ and $(\tilde q,\tilde r)$ sharp $\sigma$-admissible and $r,\tilde r < \infty$, there exists a constant $C_{2}>0$ such that for for any $\delta>0$, any $f \in L^{2}(\R^{n})$ and any $F \in L^{q'}_{t}(\R ; L^{r'}_{x}(\R^{n}))$
\be\label{estim_Strichartz}
\begin{aligned}
&\left\| e^{\iD \frac{t}{\delta} g(\delta |D|)}  \chi(\delta |D|) f \right\|_{L^{q}_{t} L^{r}_{x}} \leq C_{2}  \delta^{\beta (\frac12-\frac{1}{r})}  \| |D|^{\gamma (\frac12 - \frac{1}{r})}  \chi(\delta |D|) f \|_{L^{2}},\\
&\left\| \int_{\R} e^{-\iD \frac{s}{\delta} g(\delta |D|)} \chi(\delta |D|) F(s) ds \right\|_{L^{2}_{x}} \leq C_{2} \delta^{\beta (\frac12-\frac{1}{r})}  \| |D|^{\gamma (\frac12 - \frac{1}{r})} \chi(\delta |D|) F \|_{L^{q'}_{t} L^{r'}_{x}},\\
&\left\| \int_{\R} e^{\iD \frac{(t-s)}{\delta} g(\delta |D|)} \chi(\delta |D|) F(s) ds \right\|_{L^{q}_{t} L^{r}_{x}}  \leq C_{2}  \delta^{\beta (1-\frac{1}{r}-\frac{1}{\tilde r})} \| |D|^{\gamma (1 - \frac{1}{r} - \frac{1}{\tilde r})} \chi(\delta |D|) F \|_{L^{\tilde q'}_{t} L^{\tilde r'}_{x}}
\end{aligned}
\ee
and one also has the retarded Strichartz estimates
\be\label{estim_Strichartzretarded}
\left\| \int_{s < t} e^{\iD \frac{(t-s)}{\delta} g(\delta |D|)} \chi(\delta |D|) F(s) ds \right\|_{L^{q}_{t} L^{r}_{x}}  \leq C_{2} \delta^{\beta (1-\frac{1}{r}-\frac{1}{\tilde r})} \| |D|^{\gamma (1 - \frac{1}{r} - \frac{1}{\tilde r})} \chi(\delta |D|) F \|_{L^{\tilde q'}_{t} L^{\tilde r'}_{x}}.
\ee
Finally, one can consider the case $r=\infty$ replacing the $L^{r}$ norm by the Besov norm $\dot{B}^{0}_{\infty,2}$ and the case $\tilde r = \infty$ replacing the $L^{\tilde r'}$ norm by the Besov norm $\dot{B}^{0}_{1,2}$.
\epr
For the sake of completeness we recall that the Besov norm $\dot{B}^{0}_{r,2}$ for $r \in [1,\infty]$ is defined from a Littlewood-Paley decomposition as
\[
\left\| u \right\|_{\dot{B}^{0}_{r,2}} = \left( \sum_{j \in \mathbb{Z}} \| Q_{j}(|D|) u \|_{L^{r}}^{2} \right)^{\frac{1}{2}}.
\]

\br\label{admissiblepair}
One can also consider other ordered pairs $(p,r)$ with the restriction that $\frac{1}{q} + \frac{\sigma}{r} < \frac{\sigma}{2}$. Indeed, defining $\tilde \sigma \in (0,\sigma)$ such that $\frac{1}{q} + \frac{\tilde \sigma}{r} = \frac{\tilde \sigma}{2}$ and interpolating the assumption $I_{t,x,k}^{0}$ with the fact that $| I_{t,x,k}^{0}| \lesssim 2^{nk} \delta^{-n}$ we get
\[
\sup_{x \in \R^{n}} | I_{t,x,k}^{0} | \leq C 2^{\frac{\gamma \tilde \sigma + n(\sigma-\tilde \sigma)}{\sigma} k} \frac{\delta^{\frac{\beta \tilde \sigma + (\tilde \sigma - \sigma) n}{\sigma}}}{|t|^{\tilde \sigma}},
\]
so that one can apply the previous proposition, replacing $\sigma$ by $\tilde \sigma$, $\gamma$ by $\frac{\gamma \sigma + d(\sigma-\tilde \sigma)}{\sigma}$ and $\beta$ by $\frac{\beta \tilde \sigma + (\tilde \sigma - \sigma) n}{\sigma}$.
\er

\begin{proof}
The proof mimics the strategy used to prove Strichartz estimates for the wave equation (see for instance the Oberwolfach seminar by M. Vi\c{s}an \cite{KT_Oberwolfach} or \cite{KeelTao98}). We give the main steps. Let fix $k \in \Z_{J}$. Using the assumption on $I^{0}_{t,x,k}$ and noting that $|\mathcal{F}_{x}(I^{0}_{t,x,k})(\xi)| \lesssim 1$, we get by interpolation that for any $f \in L^{r'}(\R^{n})$
\[
\left\|  e^{\iD \frac{t}{\delta} g(\delta |D|)} P( \tfrac{\delta |D|}{2^{k}})  f \right\|_{L^{r}_{x}} \leq C_{0} 2^{\gamma (1 - \frac{2}{r})k} \frac{\delta^{\beta (1 - \frac{2}{r})}}{|t|^{\sigma (1 - \frac{2}{r})}} \| f \|_{L^{r'}}.
\]
By the Hardy-Littlewood-Sobolev inequality 
\begin{align*}
\left\| \int_{\R}  e^{\iD \frac{(t-s)}{\delta} g(\delta |D|)} P^2( \tfrac{\delta |D|}{2^{k}})F(s) ds \right\|_{L^{q}_{t} L^{r}_{x}}  &\lesssim \left\| \int \frac{\delta^{\beta (1 - \frac{2}{r})}}{|t-s|^{\sigma (1 - \frac{2}{r})}} \left\| F(s) \right\|_{L^{r'}_{x}} ds \right\|_{L^{q}_{t}}\\
&\lesssim 2^{(1 - \frac{2}{r})k} \delta^{\beta (1 - \frac{2}{r})} \left\| F(s) \right\|_{L^{q'}_{t} L^{r'}_{x}}.
\end{align*}
Then one can perform a $T^{\ast} T$ argument on the operator $T_{k}$
\[
\begin{array}{ccccc}
T_{k}& : & L^{q'}_{t}(\R, L^{r'}_{x}(\R^{n})) & \to & L^{2}(\R^{n})\\
 & & F & \mapsto & \int e^{- \iD \frac{s}{\delta} g(\delta |D|)} P(\tfrac{\delta |D|}{2^{k}})  F(s) ds
\end{array}
\]
and the first three estimates follow by applying the previous estimates to $Q_{k} f$ and $Q_{k} F$. Concerning the fourth estimate (localized retarded Strichartz estimates), we note that
\begin{align*}
\left\| \int_{s < t} e^{\iD \frac{(t-s)}{\delta} g(\delta |D|)}  Q_{k}(\delta |D|) F(s) ds \right\|_{L^{\infty}_{t} L^{2}_{x}} &= \left\| T_{k} ( \mathds{1}_{(-\infty, t)} Q_{k}(\delta|D|) F ) \right\|_{L^{\infty}_{t} L^{2}_{x}}\\
&\lesssim 2^{(\frac12-\frac{1}{\tilde r}) k}\delta^{\beta (\frac12-\frac{1}{\tilde r})}  \| Q_{k}(\delta |D|) F \|_{L^{\tilde q'}_{t} L^{\tilde r'}_{x}}.
\end{align*}
and using the strategy used to prove the first two estimates, we also have
\begin{align*}
\left\| \int_{s < t} e^{\iD \frac{(t-s)}{\delta} g(\delta |D|)}  Q_{k}(\delta |D|)  F(s) ds \right\|_{L^{\tilde q}_{t} L^{\tilde r}_{x}} &\lesssim \left\| \int \frac{\delta^{\beta (1 - \frac{2}{\tilde r})}}{|t-s|^{\sigma (1 - \frac{2}{\tilde r})}} \left\| Q_{k}(\delta |D|)  F(s) \right\|_{L^{\tilde r'}_{x}} ds \right\|_{L^{\tilde q}_{t}}\\
&\lesssim 2^{(1-\frac{2}{\tilde r}) k} \delta^{\beta (1-\frac{2}{\tilde r})}  \|  Q_{k}(\delta |D|)  F(s) \|_{L^{\tilde q'}_{t} L^{\tilde r'}_{x}}.
\end{align*}
Interpolating the two previous bounds we obtain the fourth estimate in the case $r \leq \tilde r$. Then, for $(a,b)$, $(\tilde a,\tilde b)$ sharp $\sigma$-admissible, we define the operator $S_{a,\tilde a}$
\[
\begin{array}{ccccc}
S_{k} & : & L^{b'}_{t}(\R, L^{\tilde a'}_{x}(\R^{n})) & \to & L^{b}_{t}(\R, L^{a}_{x}(\R^{n}))\\
 & & F & \mapsto & \int_{s < t} e^{\iD \frac{(t-s)}{\delta} g(\delta |D|)} P^2(\tfrac{\delta |D|}{2^{k}}) F(s) ds.
\end{array}
\]
We notice that for any $F \in L^{b'}_{t}(\R, L^{a'}_{x}(\R^{n}))$ and any function $\phi  \in L^{\tilde b'}_{t}(\R, L^{\tilde a'}_{x}(\R^{n}))$ we have
\[
\left< S_{k} F , \phi \right>_{L^{\tilde b}_{t} L^{\tilde a}_{x} \times L^{\tilde b'}_{t} L^{\tilde a'}_{x}} +  \left< F , S_{k} \phi \right>_{L^{b'}_{t} L^{a'}_{x} \times L^{b}_{t} L^{a}_{x}}  = \left< T^{\ast}_{k} T_{k} F , \phi \right>_{L^{\tilde b}_{t} L^{\tilde a}_{x} \times L^{\tilde b'}_{t} L^{\tilde a'}_{x}} 
\]
so that, using the third estimate and the previous case, the case $r > \tilde r$ follows.

For the last four estimates we sum on $\Z_{J}$ and use Bernstein's inequality and the embeddings $\dot{B}^{0}_{r,2} \hookrightarrow L^{r}$ and $L^{r'} \hookrightarrow  \dot{B}^{0}_{r',2}$ when $r \in [2,\infty[$ (see for instance Theorem 2.40 in \cite{BCD_Fourier}).
\end{proof}

The next proposition can be seen as a generalization of Theorem 2.1 in \cite{kenig_ponce_vega_kdv_disper}.

\bpr\label{Strichartz_Paleysum}
Assume that we are under the assumptions of Proposition \ref{Strichartz_Paley}, that $\sigma < 1$ and that 
\[
\left| \int_{\R^{n}} e^{\iD x \xi} e^{\iD \tfrac{t}{\delta} g(\delta |\xi|)} |\xi|^{-\gamma} \sum_{k \in \Z_{J}} Q_{k}(\delta|\xi|) d\xi \right| \leq C \frac{\delta^{\beta}}{|t|^{\sigma}}.
\]
Then one can actually take $r=\infty$ and $\tilde r=\infty$ in the bounds \eqref{estim_Strichartz} of Proposition \ref{Strichartz_Paley}. Furthermore, one can take $(r,\tilde r) = (2,\infty),(\infty,2)$ or $(\infty,\infty)$ in the retarded  Strichartz estimates \eqref{estim_Strichartzretarded} of Proposition \ref{Strichartz_Paley}.

\epr 

\br 
Together with Lemma \ref{decay1d_1deriv}, Lemma \ref{decay1d_VdC2}, Lemma \ref{decay2d_paley_1/t} and Lemma \ref{VdC2d}, we provide another way to prove, for $n=1$ or $2$, any $\alpha \notin \{ -1 , -2 \}$ and any $(q,r)$ sharp $\frac{n}{2}$-admissible, the Strichartz estimates
\begin{align*}
&\left\| e^{\iD t |D|^{\alpha+2}} |D|^{\frac{\max(\alpha,0)}{n}(\frac12-\frac{1}{r})} f \right\|_{L^{q}_{t} L^{r}_{x}} \lesssim  \| |D|^{-\frac{\min(\alpha,0)}{n}(\frac12-\frac{1}{r})} f \|_{L^{2}},\\ 
&\left\| e^{\iD t |D|^{\alpha+1} \partial_{x}} |D|^{\frac{\max(\alpha,0)}{n}(\frac12-\frac{1}{r})} f \right\|_{L^{q}_{t} L^{r}_{x}} \lesssim  \| |D|^{-\frac{\min(\alpha,0)}{n}(\frac12-\frac{1}{r})} f \|_{L^{2}} \text{   when   } n=1.
\end{align*}
\er

\br 
It is tempting to argue by interpolation as in Proposition \ref{Strichartz_Paley} and obtain more retarded Strichartz estimates. Complex interpolation requires however Banach spaces and can be tricky when one deal with $L^{1}$ or $L^{\infty}$ type function spaces, both issues we have here since we are dealing with homogeneous Sobolev spaces based on $L^{1}$. In some cases one can still interpolate using for instance arguments expounded in \cite{Gaudin_Hom_Rn} or by proving that the operators we are dealing with are part of a Stein interpolation family. Note also that such strategy provides constants $C_{1},C_{2}$ that are independent of $(q,\tilde q,r, \tilde r)$ when $\sigma<1$. Since it goes beyond the philosophy of this paper we decided not to present such estimates.
\er

\begin{proof}
We define $q_{0}>2$ such that $(q_{0},\infty)$ is sharp $\sigma$-admissible. One can use a $T^{\ast} T$ argument with the operator 
\[
\begin{array}{ccccc}
T & : & L^{q_{0}'}_{t}(\R, L^{1}_{x}(\R^{n})) & \to & L^{2}(\R^{n})\\
 & & G & \mapsto & \int \sum_{k \in \Z_{J}} Q_{k}(\delta |D|) e^{-\iD \frac{s}{\delta} g(\delta |D|)} |D|^{-\frac{\gamma}{2}} G(s) ds
\end{array}
\]
so that the first three estimates follow in the case $r=\infty$ or $\tilde r=\infty$. Concerning the fourth estimate, one can prove that
\begin{align*}
&\left\| \int_{s < t} e^{\iD \frac{(t-s)}{\delta} g(\delta |D|)} \chi(\delta |D|) F(s) ds \right\|_{L^{\infty}_{t} L^{2}_{x}} \lesssim \delta^{\frac{\beta}{2}} \left\| |D|^{\frac{\gamma}{2}} F \right\|_{L^{q_{0}'}_{t} L^{1}_{x}},\\
&\left\| \int_{s < t} e^{\iD \frac{(t-s)}{\delta} g(\delta |D|)} \chi(\delta |D|) F(s) ds \right\|_{L^{q_{0}}_{t} L^{\infty}_{x}} \lesssim \delta^{\beta} \left\| |D|^{\gamma} F \right\|_{L^{q_{0}'}_{t} L^{1}_{x}},
\end{align*}
and the bound in the case $(r,\tilde r) = (\infty,2)$ follows from a duality argument as in the proof of Proposition \ref{Strichartz_Paley}.
\end{proof}

Finally we provide Strichartz estimates that are useful for low frequencies or when $|g''| \sim 1$.

\bpr\label{Strichartz_chi}
Let $n=1$ or $2$ and $\chi$ be a smooth bounded function. Assume that $g$ satisfies (H0). Suppose there exists $C_{0}>0$, $\sigma \in (0,1]$ and $\beta, \gamma \in \R$ such that for any $\delta>0$ and any $t \in \R^{\ast}$
\[
\left| \int_{\R^{n}} e^{\iD x \xi} e^{\iD \tfrac{t}{\delta} g(\delta |\xi|)} \chi(\delta \xi) d\xi \right| \leq C_{0} \frac{\delta^{\beta}}{|t|^{\sigma}}.
\]
Then for any $(q,\tilde q,r, \tilde r) \in [2,\infty]^4$ with $(q,r)$ and $(\tilde q,\tilde r)$ sharp $\sigma$-admissible, there exists a constant $C_{2}>0$ such that for any $\delta>0$, any $f \in L^{2}(\R^{n})$ and any $F \in L^{q'}_{t}(\R ; L^{r'}_{x}(\R^{n}))$
\begin{align*}
&\left\| e^{\iD \frac{t}{\delta} g(\delta |D|)}  \chi(\delta |D|) f \right\|_{L^{q}_{t} L^{r}_{x}} \leq C_{1}  \delta^{\beta (\frac12-\frac{1}{r})}  \| f \|_{L^{2}},\\
&\left\| \int_{\R} e^{-\iD \frac{s}{\delta} g(\delta |D|)}  \chi(\delta |D|) F(s) ds \right\|_{L^{2}_{x}} \leq C_{1} \delta^{\beta (\frac12-\frac{1}{r})}  \| F \|_{L^{q'}_{t} L^{r'}_{x}},\\
&\left\| \int e^{\iD \frac{(t-s)}{\delta} g(\delta |D|)} \chi^2(\delta |D|)  F(s) ds \right\|_{L^{q}_{t} L^{r}_{x}}  \leq C_{2}  \delta^{\beta (1-\frac{1}{r}-\frac{1}{\tilde r})} \| F \|_{L^{\tilde q'}_{t} L^{\tilde r'}_{x}},\\
&\left\| \int_{s < t} e^{\iD \frac{(t-s)}{\delta} g(\delta |D|)}  \chi^2(\delta |D|)  F(s) ds \right\|_{L^{q}_{t} L^{r}_{x}}  \leq C_{2} \delta^{\beta (1-\frac{1}{r}-\frac{1}{\tilde r})} \| F \|_{L^{\tilde q'}_{t} L^{\tilde r'}_{x}}.
\end{align*}
\epr

\begin{proof}
The proposition follows from a $T^{\ast} T$ argument and a duality argument on the operator 
\[
\begin{array}{ccccc}
T & : & L^{q'}_{t}(\R, L^{r'}_{x}(\R^{n})) & \to & L^{2}(\R^{n})\\
 & & F & \mapsto & \int e^{- \iD \frac{s}{\delta} g(\delta |D|)} \chi(\delta |D|)  F(s) ds.
\end{array}
\]
\end{proof}

\subsection{Local Kato smoothing/Morawetz type estimates}\label{ss:kato-morawetz}

Some dispersive propagators provide a local smoothing effect (also called local Kato smoothing effect \cite{Kato83}) or some uniform bounds of the local energy (also called Morawetz estimates \cite{Morawetz68}). The next proposition gives a unified version of these two properties.

\bpr\label{kato-morawetz}
Let $n \in \N^{\ast}$. Assume that $g$ is $\mathcal{C}^{1}(\R^{\ast}_{+})$. There exists $C>0$ such that for any function $f$ in $L^{2}(\R^{n})$, any $\delta>0$, any $a>0$ and any $x_{0} \in \R^{n}$
\[
\int_{\R} \int_{\R^{n}} \left| (|g'(\delta |D|)|^{\frac12} e^{\iD \frac{t}{\delta} g(\delta |D|)} f)(x) \right|^2 e^{-\frac{a}{2} |x-x_{0}|^2} dx dt \leq \frac{C}{\sqrt{a}} \| f \|_{L^{2}}^2.
\]
\epr

\br 
This proposition can be read as follows. Assume by simplicity that $\delta=1$.
\begin{enumerate}
\item If $|g'(y)| \sim y^{\beta}$ on $\R^{+}$ with $\beta>0$, we obtain a local Kato smoothing
\[
\int_{\R} \int_{\R^{n}} \left| (|D|^{\frac{\beta}{2}} e^{\iD t g(|D|)} f)(x) \right|^2 e^{- |x|^2} dx dt \lesssim \| f \|_{L^{2}}^2.
\]
$\beta=2$ (KdV-type behavior) corresponds to \cite{Kato83}, $\beta=1$ (Schr\"odinger-type behavior) to \cite{ConstantinSaut88} and $\beta=\frac{1}{2}$ (water-wave with surface tension and infinite depth) to \cite{ABZ10}. 
\item If $|g'(y)| \sim 1$ on $\R^{+}$ (wave-type behavior), we get a Morawetz type estimate  
\[
\sup_{x_{0} \in \R^{n}} \int_{\R} \int_{\R^{n}} \left| (e^{\iD t g(|D|)} f )(x)\right|^2 e^{- |x-x_{0}|^2} dx dt \lesssim \| f \|_{L^{2}}^2.
\]
\item If $|g'(y)| \sim y^{\beta}$ on $\R^{+}$ with $\beta<0$ (weakly dispersive phase), we have a Morawetz type estimate
\[
\sup_{x_{0} \in \R^{n}} \int_{\R} \int_{\R^{n}} \left| (e^{\iD t g(|D|)} f)(x) \right|^2 e^{- |x-x_{0}|^2} dx dt \lesssim \| |D|^{-\frac{\beta}{2}} f \|_{L^{2}}^2.
\]
Such a kind of estimate was obtained for instance with $\beta=-\frac{1}{2}$ (water-wave with infinite depth) in \cite{AIT22}.

Finally, note that the bound is uniform with respect to $\delta$ so that if $g'(0) \neq 0$, the case $\delta=0$ exactly corresponds to a Morawetz-type estimate for the wave propagator.
\end{enumerate}
\er

We only chose to present one type of Morawetz estimates (with a Gaussian weight). Other types of Morawetz estimates can be obtained adapting for instance \cite{OR2013} to nonhomogeneous radial phases. 

\begin{proof}
The strategy used here is similar to the proof of the local smoothing of the Schr\"odinger propagator in the Oberwolfach seminar by M. Vi\c{s}an \cite{KT_Oberwolfach} (see Lemma 2.11). By a change of variables on can assume that $\delta=1$. We argue by density and assume that $f$ is a Schwartz class function whose Fourier transform support does not contain $0$. We notice that
\[
\int_{\R} \int_{\R^{n}} \left| (e^{\iD t g(|D|)}  |g'(|D|)|^{\frac12} f)(x) \right|^2 e^{- \frac{a}{2} |x-x_{0}|^2} dx dt = \int_{t \in \R} \hspace{-0.05cm} \int_{\xi \in \R^{n}} \hspace{-0.05cm} \int_{\eta \in \R^{n}} \hspace{-0.2cm} G_{a}(\eta,\xi) e^{\iD t (g(|\xi|)-g(|\eta|))} d\eta d\xi dt
\]
where
\[
G_{a}(\eta,\xi) := \frac{1}{(2 \pi a)^{\frac{n}{2}}} e^{-\iD x_{0} \cdot (\xi-\eta)} e^{-\frac{1}{2 a} |\xi-\eta|^{2}} |g'(|\xi|)|^{\frac12} |g'(|\eta|)|^{\frac12} \widehat{f}(\xi) \overline{\widehat{f}(\eta)}.
\]
Let $J \subset \R^{+}$ be an interval such that $|g'|>0$ on $J$. We define
\[
I_{J} :=  \int_{t \in \R} \int_{\xi \in \R^{n}} \int_{\eta \in \R^{n}} G_{a}(\eta,\xi) \mathds{1}_{\{ |\xi| \in J\}} \mathds{1}_{\{ |\eta| \in J\}} e^{\iD t (g(|\xi|)-g(|\eta|))} d\eta d\xi dt.
\]
We get
\begin{align*}
I_{J} &= \int_{\xi \in \R^{n}} \int_{\eta \in \R^{n}} G_{a}(\eta,\xi) \mathds{1}_{\{ |\xi| \in J\}} \mathds{1}_{\{ |\eta| \in J\}} \int_{t \in \R}  e^{\iD t (g(|\xi|)-g(|\eta|))} dt d\eta d\xi\\
&= \sqrt{2 \pi} \int_{\xi \in \R^{n}} \int_{\eta \in \R^{n}} \delta_{\{ g(|\xi|) = g(|\eta|)\}} G_{a}(\eta,\xi) \mathds{1}_{\{ |\xi| \in J\}} \mathds{1}_{\{ |\eta| \in J\}} d\eta d\xi\\
&= \sqrt{2 \pi}\int_{\xi \in \R^{n}} \int_{\eta \in \R^{n}} \delta_{\{|\xi| = |\eta| \}} \frac{G_{a}(\eta,\xi)}{|g'(|\xi|)|} \mathds{1}_{\{ |\xi| \in J\}} \mathds{1}_{\{ |\eta| \in J\}} d\eta d\xi
\end{align*}
where $\delta_{\{|\xi| = |\eta| \}}$ is the Dirac delta function. Then, using polar coordinates $\xi=r \omega$ and $\eta=y u$ and denoting by $S^{n-1}$ the $n-1$ sphere
\begin{align*}
I_{J} &= \sqrt{2 \pi} \int_{r=0}^{\infty} \int_{y=0}^{\infty} \int_{\omega \in S^{n-1}} \int_{u \in S^{n-1}} \delta_{\{r = y \}} \frac{G_{a}(y u,r \omega) }{|g'(r)|} \mathds{1}_{\{ r \in J\}} \mathds{1}_{\{y \in J\}} du d\omega dy dr\\
&= \frac{(2 \pi)^{\frac{1-n}{2}}}{a^{\frac{n}{2}}} \int_{r=0}^{\infty} \int_{\omega \in S^{n-1}} \int_{u \in S^{n-1}} e^{-\iD r x_{0} \cdot (\omega-u)}  e^{-\frac{r^2}{2 a} |\omega-u|^{2}} \mathds{1}_{\{ r \in J\}} \widehat{f}(r \omega) \overline{\widehat{f}(y u)} r^{2(n-1)}du d\omega dr.
\end{align*}
We finally notice that if we define $K_{a}(r,r\omega,ru) := \frac{1}{a^{\frac{n}{2}}} e^{-\iD r x_{0} \cdot (\omega-u)} e^{-\frac{r^2}{2 a} |\omega-u|^{2}} r^{n-1}$
\[
\left| \int_{\omega \in S^{n-1}} K_{a}(r,r\omega,ru) d\omega \right| + \left| \int_{u \in S^{n-1}} K_{a}(r,r\omega,ru) du \right| \lesssim \frac{1}{a^{\frac{n}{2}}} \int_{0}^{\pi} e^{- \tfrac{1}{5 a} r^{2} \theta^{2}} r^{n-1} \theta^{n-2} d\theta \lesssim \frac{1}{\sqrt{a}}.
\]
The result follows from Schur's test.
\end{proof}

We can now state local energy decay.

\bc\label{decay_local_energy}
Let $n \in \N^{\ast}$. Assume that $g$ is $\mathcal{C}^{1}(\R^{\ast}_{+})$. Let $a>0$, $\delta>0$, $x_{0} \in \R^{n}$ and $f \in L^{2}(\R^{n})$. Then
\[
\lim_{t \to \pm \infty} \int_{\R^{n}} \left| (|g'(\delta |D|)|^{\frac12} e^{\iD \frac{t}{\delta} g(\delta |D|)} f)(x) \right|^2 e^{-\frac{a}{2} |x-x_{0}|^2} dx = 0.
\]
Furthermore, if $g$ is defined and $\mathcal{C}^1$ in the vicinity of $0$, the limit is uniform with respect to $\delta \to 0$.
\ec

\begin{proof}
We argue by density and assume that $f$ is a Schwartz class function whose Fourier transform support does not contain $0$ and is bounded. We define the map
\[
E_{a,x_{0},\delta} : t \mapsto \int_{\R^{n}} \left| (|g'(\delta |D|)|^{\frac12} e^{\iD \frac{t}{\delta} g(\delta |D|)} f)(x) \right|^2 e^{-\frac{a}{2} |x-x_{0}|^2} dx. 
\]
Using Proposition \ref{kato-morawetz} and since $\| f \|_{L^{2}}$ and $\| \tfrac{g(\delta |D|)}{\delta} f \|_{L^{2}}$ are finite, the maps $E_{a,x_{0},\delta}$ and $E_{a,x_{0},\delta}'$ are $L^{1}(\R)$ so that $E_{a,x_{0},\delta}$ goes to $0$ at $\pm \infty$. Finally, when $g$ is defined and $\mathcal{C}^1$ in the vicinity of $0$, $\| f \|_{L^{2}}$ and $\| \tfrac{g(\delta |D|)-g(0)}{\delta} f \|_{L^{2}}$ are bounded uniformly with respect to $\delta \to 0$.
\end{proof}

Actually when $n=1$ one can prove a stronger version of Proposition \ref{kato-morawetz}.

\bpr\label{kato-morawetz_n=1}
Let $n=1$. Assume that $g$ is $\mathcal{C}^{1}(\R^{\ast}_{+})$. Then, there exists a constant $C>0$ such that for any function $f$ in $L^{2}(\R^{n})$ and any $\delta>0$
\[
\sup_{x \in \R} \int_{\R} \left| (| g'(\delta |D|)|^{\frac12} e^{\iD \frac{t}{\delta} g(\delta |D|)} f)(x) \right|^2 dt \leq C \| f \|_{L^{2}}^2.
\]
\epr

This proposition can be seen as a generalization of \cite{Vega88} and a reformulation of Theorem 4.1 in \cite{kenig_ponce_vega_kdv_disper}.

\begin{proof}
For the sake of completeness we recall the proof. Firstly using a change of variables one can assume $\delta=1$. Then, let $J$ be an interval such that $|g'|>0$ on $J$. Therefore $g$ is invertible on $J$ and by simplicity, we denote in the following as $g^{-1}$ the inverse of $g$ on $J$. We define
\begin{align*}
u_{J}(t,x) &:= \frac{1}{\sqrt{2\pi}} \int_{\R} |g'( |\xi|)|^{\frac12} e^{\iD t g(|\xi|)} e^{\iD x \xi} \widehat{f}(\xi) \mathds{1}_{\{ \xi \in J \} } d\xi\\
&= \frac{1}{\sqrt{2\pi}} \int_{\R}  e^{\iD t y} |g'( |g^{-1}(y)|)|^{-\frac12} e^{\iD x g^{-1}(y)} \widehat{f}(g^{-1}(y)) \mathds{1}_{\{ y \in g(J) \} }  dy.
\end{align*}
The key observation is to see $u$ as a Fourier transform with respect to the variable $t$ so that using Plancherel's identity
\[
\int_{\R} |u_{J}(t,x)|^2 dt \lesssim \int_{\R}  |g'( |g^{-1}(y)|)|^{-1} |\widehat{f}(g^{-1}(y))|^2 \mathds{1}_{\{ y \in g(J) \} } dy = \int_{\R}  |\widehat{f}(\xi)|^2 \mathds{1}_{\{ \xi \in J \} } d\xi.
\]
\end{proof}

\section{Examples}
In this section we apply the results obtained in the previous section to various dispersive equations. We recall that we study equations under the form
\bes
\partial_{t} u = \pm \frac{\iD}{\delta} g(\delta |D|) u,
\ees
or when $n=1$
\bes
\partial_{t} u = \pm \frac{g(\delta |D|)}{\delta |D|} \partial_{x} u.
\ees
We use the notations of Section \ref{sec_dispersive}. Notice that for any smooth bounded function $\tilde \chi$, $s \in \R$, $n=1$ or $2$, if we define
\[
A = \{ k \in \Z \text{ , } \exists  r \in [\tfrac12,2] \text{ such that } \tilde \chi(2^{k} r) \neq 0 \}
\]
and if we take $P=Q_{0}$ in the definition of $I_{t,x,k}$, then 
\[
\int_{\R} e^{\iD x \cdot \xi} e^{\iD \tfrac{t}{\delta} g(\delta |\xi|)} \tilde \chi(\delta |\xi|) \widehat{u_{0}}(\xi) d\xi = (2\pi)^{\frac{n}{2}} \sum_{k \in A} (I^{s}_{t,\cdot,k} -  I^{s}_{-t,- \cdot,k}) \ast |D|^{-s} \tilde \chi(\delta |D|)  u_{0}.
\]

\subsection{Ostrovsky equation}\label{ss:ostrovsky}

We consider the linear Ostrovsky equation
\[
\partial_{t} u = ( - \partial_{x}^{-1} u + b \partial_{x}^{3} u )
\]
so that $g(y) =\frac{1}{y} - b y^{3}$ with $b \in \R^{\ast}$ and $\delta=1$. Dispersive estimates were obtained in \cite{varlamov_ostrov} when $b>0$ and some Strichartz estimates and a local Kato smoothing were established in \cite{linares_ostrov} for any $b \in \R^{\ast}$.

If $b<0$, $g''>0$ on $\R^{\ast}_{+}$ and Lemma \ref{decay1d_1deriv}, Lemma \ref{decay1d_VdC2} and Proposition \ref{Strichartz_Paleysum} give for any $t \in \R^{\ast}$ and any $(q,r)$ sharp $\tilde \sigma$-admissible with $\tilde \sigma \in [\frac13,\frac12]$
\[
\left\| e^{\iD t g(D)} u_{0} \right\|_{L^{\infty}_{x}} \lesssim \frac{\| u_{0}\|_{L^{1}}}{\sqrt{|t|}} \text{   ,   } \left\| e^{\iD t g(D)} u_{0} \right\|_{L^{q}_{t} L^{r}_{x}} \lesssim \| u_{0} \|_{L^{2}}
\]  
and we can obtain a local Kato smoothing thanks to Proposition \ref{kato-morawetz_n=1} as in \cite{linares_ostrov}.

If $b>0$, there exists a unique $y_{b}>0$ such that $g''(y_{b})=0$ and $g'' \leq -6b$ on $\R^{+}$. Introducing a smooth compactly supported function $\tilde \chi$ that is equal to $1$ near $\pm y_{b}$ and whose support does not contain $0$, Lemma \ref{decay1d_VdC} gives for any $t \in \R^{\ast}$
\[
\left\| e^{\iD t g(D)} \tilde \chi(|D|) u_{0} \right\|_{L^{\infty}_{x}} \lesssim \frac{\| u_{0}\|_{L^{1}}}{|t|^\frac13} \text{  ,  } \left\| e^{\iD t g(D)} (1-\tilde \chi(|D|)) u_{0} \right\|_{L^{\infty}_{x}} \lesssim \min \left( \frac{1}{|t|^{\frac13}} , \frac{1}{\sqrt{|t|}} \right) \| u_{0} \|_{L^{1}}
\]
and for any $(q,r)$ sharp $\frac12$-admissible and any $(\tilde q, \tilde r)$ sharp $\frac13$-admissible, Proposition \ref{Strichartz_chi} provides the Strichartz estimates
\[
\left\| e^{\iD t g(D)} (1-\tilde \chi(|D|)) u_{0} \right\|_{L^{q}_{t} L^{r}_{x}} \lesssim \| u_{0} \|_{L^{2}} \text{   and   } \left\| e^{\iD t g(D)} u_{0} \right\|_{L^{\tilde q}_{t} L^{\tilde r}_{x}} \lesssim \| u_{0} \|_{L^{2}}.
\]
One can also prove a local Kato smoothing thanks to Proposition \ref{kato-morawetz_n=1} as in \cite{linares_ostrov}.

Finally one can deal with the case $b=0$. It is a linearized version of an equation which has several names in the literature : Vakhnenko equation (\cite{Vakhnenko92}), Ostrovsky-Hunter equation (\cite{Hunter90}) or simply reduced Ostrovsky equation. We get a $L^{\infty}$ dispersive estimate from Lemma \ref{decay1d_VdC2}, Strichartz estimates from Proposition \ref{Strichartz_Paleysum}  and Morawetz type estimates by Proposition \ref{kato-morawetz_n=1}: for any $t \in \R^{\ast}$ and any $(q,r)$ sharp $\tfrac12$-admissible
\begin{align*}
&\left\| e^{\frac{\iD t}{|D|}} u_{0} \right\|_{L^{\infty}_{x}} \lesssim \frac{\| |D|^{\frac32} u_{0}\|_{L^{1}}}{\sqrt{|t|}} \text{  ,  } \left\| e^{\frac{\iD t}{|D|}} u_{0} \right\|_{L^{q}_{t} L^{r}_{x}} \lesssim \| |D|^{\frac{3}{2}(\frac{1}{2}-\frac{1}{r})} u_{0} \|_{L^{2}},\\
&\sup_{x \in \R} \int_{\R} |(e^{\frac{\iD t}{|D|}} \tfrac{1}{|D|} u_{0})(x)|^2 dt \lesssim \| u_{0} \|_{L^{2}} \text{    ,    } \sup_{x \in \R} \int_{\R} |(e^{\frac{\iD t}{|D|}} u_{0})(x)|^2 dt \lesssim \| |D| u_{0} \|_{L^{2}}.
\end{align*}

\subsection{BBM-KdV equation}\label{ss:BBM_KdV}

We consider a linear BBM-KdV equation (named after \cite{BBM})
\[
\partial_{t} u + \mu p \partial_{x}^{2} \partial_{t} u = \pm (\partial_{x} u + \mu (p+\tfrac16) \partial_{x}^{3} u).
\]
Here $g(y) = y \frac{1 - (p+\tfrac16) y^2}{1- p y^2}$ with $p \leq 0$ and $\delta=\sqrt{\mu} \in (0,1]$. The equation we present here is the linear version of (7.7) in \cite{Lannes_ww} and serves as a model for the propagation of long waves. Since $p=0$ corresponds to the KdV case that is well understood we focus on the case $p<0$. We note that
\[
g'(y) = \frac{6-(12p+3)y^2+p(6p+1)y^4}{6(1-py^2)^2}\text{  ,  } g''(y) = - \frac{y(3+p y^2)}{3(1-py^{2})^3} \text{  ,  } g'''(y)= - \frac{p^2 y^4 + 6p y^2 + 1}{(1-py^2)^4}
\]
so that $g''$ has a unique positive zero $y_{0} :=\sqrt{\tfrac{3}{|p|}}$, $g''$ and $g'''$ have no common zero and
\[
g'(y)-\frac{6p+1}{6p} \underset{+\infty}{\sim} -\frac{1}{6 p^2 y^{2}} \text{   ,   } g''(y) \underset{+\infty}{\sim} \frac{1}{3 p^2 y^{3}}.
\]
We introduce a smooth compactly supported function $\chi$ that is equal to $1$ on $[-3y_{0},3y_{0}]$. From Lemma \ref{decay1d_VdC} (for the low frequencies), Lemma \ref{decay1d_multideriv} (for the intermediate frequencies), Lemma \ref{decay1d_VdC2} (for the high frequencies) and Proposition \ref{Strichartz_Paleysum}, for any $l \geq 2$, any $(\tilde q, \tilde r)$ sharp $\frac12$-admissible, any $(\tilde q, \tilde r)$ sharp $\frac13$-admissible, there exists $C>0$ such that  for any $t \in \R^{\ast}$ and any $\mu \in (0,1]$
\begin{align*}
&\left\| e^{\iD \tfrac{t}{\sqrt{\mu}} g(\sqrt{\mu} D)} \chi(\sqrt{\mu} D) u_{0} \right\|_{L^{\infty}_{x}} \leq C \frac{\| u_{0}\|_{L^{1}}}{(\mu |t|)^{\frac13}},\\
&\left\| e^{\iD \tfrac{t}{\sqrt{\mu}} g(\sqrt{\mu} D)} (1-\chi(\sqrt{\mu} D))  u_{0} \right\|_{L^{\infty}_{x}} \leq C \frac{\mu^{\frac{1}{l}}}{|t|^{\frac{1}{l}}} \| |D|^{\frac{1+l}{l}} (1-\chi(\sqrt{\mu} D)) u_{0} \|_{L^{1}},\\
&\left\| e^{\iD \tfrac{t}{\sqrt{\mu}} g(\sqrt{\mu} D)} (1-\chi(\sqrt{\mu} D))  u_{0} \right\|_{L^{q}_{t} L^{r}_{x}} \leq C \mu^{\frac{1}{2}(\frac12-\frac{1}{r})} \| |D|^{\frac32 (\frac12-\frac{1}{r})} (1-\chi(\sqrt{\mu} D)) u_{0} \|_{L^{2}},\\
&\left\| e^{\iD \tfrac{t}{\sqrt{\mu}} g(\sqrt{\mu} D)} u_{0} \right\|_{L^{\tilde q}_{t} L^{\tilde r}_{x}} \leq \frac{C}{\mu^{\frac13 (\frac12 - \frac{1}{\tilde r})}} \| u_{0} \|_{L^{2}} + C \mu^{\frac13 (\frac12 - \frac{1}{\tilde r})} \| (1-\chi(\sqrt{\mu} D)) |D|^{\frac43 (\frac12 - \frac{1}{\tilde r})} u_{0} \|_{L^{2}}.
\end{align*}
Finally one can prove a Morawetz type estimate for some $p$. We notice that $g'$ has a positive root if and only if $p \geq -\tfrac{3}{16}$. Then, if $p<-\tfrac{3}{16}$, from Proposition \ref{kato-morawetz_n=1} and Corollary \ref{decay_local_energy}, there exists $C>0$ such that for any $u_{0} \in L^{2}(\R^{n})$ and any $\mu \in (0,1]$
\[
\sup_{x \in \R} \int_{\R} |(e^{\iD \frac{t}{\sqrt{\mu}} g(\sqrt{\mu} |D|)} u_{0})(x)|^2 dt \leq C \left\| u_{0} \right\|_{L^{2}}^2 \text{   ,   } \lim_{t \to \pm \infty} \int_{\R} |(e^{\iD \frac{t}{\sqrt{\mu}} g(\sqrt{\mu} |D|)} u_{0})(x)|^2 e^{-|x|^2} dx = 0.
\]

\subsection{Intermediate long wave equation}\label{ss:ILW}

We consider the linear intermediate long wave equation
\[
u_{t} =  \frac{1}{\rho} \varphi(\rho |D|) \partial_{x} u
\]
where $\varphi(y) = y \coth(y) - 1$ and with $\rho>0$. We refer for instance to \cite{book_Klein_Saut} for more details concerning the nonlinear version of this equation. Note that one must change $t$ by $\tfrac{t}{\rho}$ in order to be under the form of Remark \ref{other_eq_n=1} with $\delta=\rho$. We define $g(y) := y \varphi(y)$ and we notice that
\begin{align*}
&g'(y)=y^{2}(1-\coth^2(y))+2y\coth(y)-1,\\
&g''(y)=2 \coth(y)+4y(1-\coth^2(y))-2y^{2} \coth(y) (1-\coth^2(y)),\\
&g''(y) \underset{0}{\sim} 2y \text{   ,   } g''(y) \underset{+\infty}{\sim} 2 \text{   ,   } g'(y) \underset{0}{\sim} y^2 \text{   ,   } g'(y) \underset{+\infty}{\sim} 2y.
\end{align*}
Note also that the function $g''$ is positive on $\R^{\ast}_{+}$ and $g'''>0$ on $[-2,2]$. We introduce a smooth compactly supported function $\chi$ that is equal to $1$ on $[-1,1]$ and whose support is a subset of $[-2,2]$. We get from Lemma \ref{decay1d_VdC},  Lemma \ref{decay1d_VdC2} and Proposition \ref{Strichartz_Paleysum} that for any $l \geq 2$, any $(q,r)$ sharp $\frac12$-admissible and any $(\tilde q, \tilde r)$ sharp $\frac13$-admissible, there exists $C>0$ such that for any $t \in \R^{\ast}$ and any $\rho>0$
\begin{align*}
&\left\| e^{\frac{t}{\rho} \varphi(\rho D) \partial_{x}} \chi(\rho D) u_{0} \right\|_{L^{\infty}_{x}} \leq C \frac{\| \chi(\rho D)  u_{0}\|_{L^{1}}}{(\rho |t|)^{\frac13}},\\
&\left\| e^{\frac{t}{\rho} \varphi(\rho D) \partial_{x}} (1-\chi(\rho D)) u_{0} \right\|_{L^{\infty}_{x}} \leq C \frac{1}{|t|^{\frac{1}{l}}} \| (1-\chi(\rho D))  |D|^{\frac{l-2}{l}} u_{0}\|_{L^{1}},\\
&\left\| e^{\frac{t}{\rho} \varphi(\rho D) \partial_{x}} (1-\chi(\rho D)) u_{0} \right\|_{L^{q}_{t} L^{r}_{x}} \leq C \| (1-\chi(\rho D)) u_{0}\|_{L^{2}},\\
&\left\| e^{\frac{t}{\rho} \varphi(\rho D) \partial_{x}} u_{0} \right\|_{L^{\tilde q}_{t} L^{\tilde r}_{x}} \leq \frac{C}{\rho^{\frac13 (\frac12 - \frac{1}{\tilde r})}} \| \chi(\rho D)  u_{0}\|_{L^{2}} + C \| (1-\chi(\rho D))  |D|^{\frac13 (\frac12 - \frac{1}{\tilde r})} u_{0}\|_{L^{2}}.
\end{align*}
We can also prove a local smoothing. Since $g'>0$, by Proposition \ref{kato-morawetz_n=1} there exists $C>0$ such that for any $u_{0} \in L^{2}(\R^{n})$ and any $\rho>0$
\[
\sup_{x \in \R} \int_{\R}  |(e^{\iD \frac{t}{\rho} \varphi(\rho |D|) \partial_{x}} ( \mathds{1}_{\{ \rho |D|<1 \}} \sqrt{\rho} |D|+\mathds{1}_{\{ \rho |D|>1 \}} \sqrt{|D|} ) u_{0})(x)|^2 dt \leq C \| u_{0} \|_{L^{2}}^2.
\]
Note that by letting $\rho \to \infty$ one get the dispersive estimates of the Benjamin-Ono equation $u_{t} = |D| \partial_{x} u$ (see for instance \cite{book_Klein_Saut}).

\subsection{Irrotational water wave equations}\label{ss:ww}

We consider the equation
\[
\partial_t u  = \pm \iD \sqrt{\frac{\tanh(\sqrt{\mu} |D|)}{\sqrt{\mu} |D|}} |D| u.
\]
Here $g(y) = \sqrt{y \tanh(y)}$, $\delta=\sqrt{\mu} \in (0,1]$ and $n=1$ or $2$. This equation arises as a factorized form of System \eqref{linear_ww}. Some $L^{1} \to L^{\infty}$ weighted estimates of the propagator were performed in \cite{my_proud_res,mesognon_dispersion}, localized dispersive and Strichartz estimates when $\mu=1$ were obtained in \cite{wang_globalWW,dinvay_et_al20,deneke_et_al22} and a Morawetz-type bound when $\mu=1$ and in the infinite depth limit was proved in \cite{AIT22}. 

We note that $g'>0$, $g''<0$ and
\[
g'(y) -1 \underset{0}{\sim} - \tfrac12 y^2 \text{  ,  } g'(y) \underset{+\infty}{\sim} \tfrac12 y^{-\tfrac12} \text{  ,  } g''(y) \underset{0}{\sim} - y \text{  ,  } g''(y) \underset{+\infty}{\sim} -\tfrac14 y^{-\tfrac32} \text{  ,  } g'''(y) \underset{+\infty}{\sim} \tfrac38 y^{-\tfrac52}.
\]
Let introduce a smooth compactly supported function $\chi$  that is equal to $1$ near $0$.

When $n=1$, we get from Lemma \ref{decay1d_VdC2} and Proposition \ref{Strichartz_Paleysum} that for any $l \geq 2$, any $(q,r)$ sharp $\frac12$-admissible and any $(\tilde q, \tilde r)$ sharp $\frac13$-admissible, there exists $C>0$ such that for any $t \in \R^{\ast}$ and any $\mu \in (0,1]$
\begin{align*}
&\left\|e^{\iD \tfrac{t}{\sqrt{\mu}} g(\sqrt{\mu} |D|)} \chi(\sqrt{\mu} D) u_{0} \right\|_{L^{\infty}_{x}} \leq C \frac{\|u_{0}\|_{L^{1}}}{(\mu |t|)^{\frac13}},\\
&\left\|e^{\iD \tfrac{t}{\sqrt{\mu}} g(\sqrt{\mu} |D|)} (1-\chi(\sqrt{\mu} D)) u_{0} \right\|_{L^{\infty}_{x}}  \leq C \frac{\mu^{\frac{1}{4l}}}{|t|^{\frac{1}{l}}} \| |D|^{\frac{2l-1}{2l}} (1- \chi(\sqrt{\mu} D)) u_{0} \|_{L^{1}},\\
&\left\|e^{\iD \tfrac{t}{\sqrt{\mu}} g(\sqrt{\mu} |D|)}  (1 - \chi(\sqrt{\mu} D)) u_{0} \right\|_{L^{q}_{t} L^{r}_{x}} \leq C \mu^{\frac18(\frac12-\frac1r)} \| |D|^{\frac34(\frac12-\frac{1}{r})} (1 - \chi(\sqrt{\mu} D)) u_{0} \|_{L^{2}},\\
&\left\|e^{\iD \tfrac{t}{\sqrt{\mu}} g(\sqrt{\mu} |D|)} u_{0} \right\|_{L^{\tilde q}_{t} L^{\tilde r}_{x}} \leq \frac{C}{\mu^{\frac13 (\frac12-\frac{1}{\tilde r})}} \| \chi(\rho D)  u_{0}\|_{L^{2}} + \mu^{\frac{1}{12}(\frac12-\frac{1}{\tilde r})} \| |D|^{\frac56(\frac12-\frac{1}{\tilde r})} (1 - \chi(\sqrt{\mu} D))  u_{0} \|_{L^{2}}.
\end{align*}
Retarded Strichartz estimates can also be obtained.

When $n=2$, we get from Lemma \ref{decay2d_ww_lowfreq}, Lemma \ref{VdC2d}, Proposition \ref{Strichartz_Paleysum}, Lemma \ref{decay2d_wave_2deriv} and Proposition \ref{Strichartz_Paley} that for any $(q,r)$ sharp $1$-admissible and any $(\tilde q, \tilde r)$ sharp $\frac12$-admissible with $\tilde r < \infty$, there exists $C>0$ such that for any $t \in \R^{\ast}$ and any $\mu \in (0,1]$
\begin{align*}
&\left\|e^{\iD \tfrac{t}{\sqrt{\mu}} g(\sqrt{\mu} |D|)} u_{0} \right\|_{L^{\infty}_{x}} \hspace{-0.3cm} \leq \frac{C}{\sqrt{\mu} |t|} \| \chi(\sqrt{\mu}  |D|) u_{0}\|_{L^{1}} \hspace{-0.05cm}  + \hspace{-0.05cm} C \frac{\mu^{\frac14}}{|t|}\| |D|^{\frac32} (1- \chi(\sqrt{\mu}|D|)) u_{0}\|_{L^{1}},\\
&\left\|e^{\iD \tfrac{t}{\sqrt{\mu}} g(\sqrt{\mu} |D|)} u_{0} \right\|_{L^{q}_{t} L^{r}_{x}} \hspace{-0.5cm} \leq \frac{C}{\mu^{\frac12(\frac12 - \frac1r)}} \| \chi(\sqrt{\mu}  |D|) u_{0} \|_{L^{2}} \hspace{-0.05cm} + \hspace{-0.05cm} C \mu^{\frac14(\frac12 - \frac1r)} \| |D|^{\frac32 (\frac12-\frac1r)} (1- \chi(\sqrt{\mu}|D|)) u_{0}\|_{L^{2}}, \\
&\left\|e^{\iD \tfrac{t}{\sqrt{\mu}} g(\sqrt{\mu} |D|)} u_{0} \right\|_{L^{\tilde q}_{t} L^{\tilde r}_{x}} \hspace{-0.5cm} \leq C \| |D|^{\frac32 (\frac12 - \frac{1}{\tilde r})}  \chi(\sqrt{\mu}  |D|) u_{0} \|_{L^{2}} \hspace{-0.05cm} + \hspace{-0.05cm} C \mu^{\frac18 (\frac12 - \frac{1}{\tilde r})} \| |D|^{\frac74 (\frac12 - \frac{1}{\tilde r})} (1- \chi(\sqrt{\mu}|D|)) u_{0} \|_{L^{2}}.
\end{align*}
Retarded Strichartz estimates can also be obtained.

Finally, one can prove a Morawetz type estimate. Since $g'>0$, from Proposition \ref{kato-morawetz} and Corollary \ref{decay_local_energy} there exists $C>0$ such that for any $u_{0} \in H^{\frac12}(\R^{n})$ and any $\mu \in (0,1]$
\begin{align*}
&\sup_{x_{0} \in \R^{n}} \int_{\R} \int_{\R^{n} } |(e^{\iD \frac{t}{\sqrt{\mu}} g(\sqrt{\mu} |D|)} u_{0})(x)|^2 e^{-|x-x_{0}|^{2}} dx dt \leq C \left\| \sqrt{1+\sqrt{\mu} |D|} u_{0} \right\|_{L^{2}}^2,\\
&\lim_{t \to \pm \infty} \int_{\R^{n}} \left| (e^{\iD \frac{t}{\sqrt{\mu}} g(\sqrt{\mu} |D|)} u_{0})(x) \right|^2 e^{- |x|^2} dx = 0.
\end{align*}

\subsection{abcd-Boussinesq systems}\label{ss:abcd_Boussi}

We consider the equation
\[
\partial_t u  = \pm \iD \varphi(\sqrt{\mu} |D|) |D| u,
\]
where 
\[
\varphi(y) = \sqrt{\frac{(1-\mu a y^2)(1-\mu c y^2)}{(1+\mu b y^2)(1+\mu d y^2)}}
\]
and here $g(y) = y \varphi(y)$, $\delta=\sqrt{\mu}$ and $n=1$ or $2$. This equation arises as a factorized form of System \eqref{linear_abcd}. We assume that
\be\label{condition_abcd}
b \geq 0 \text{  ,  } d \geq 0 \text{  ,  } a \leq 0 \text{  ,  } c \leq 0
\ee
in order to get a wellposed equation (see \cite{bona_chen_saut_derivation}) and that
\be\label{condition_abcd2}
((a+b)(a+d)(c+b)(c+d))^2+(a+b+c+d)^2 > 0
\ee
so that we avoid the situation where $g(r) \equiv r $ which corresponds to the half-wave case. Some dispersive estimates were obtained in the case $a=c=\frac16$ and $b=d=0$ in \cite{LPS12}.

We note that
\bes
g'(y) = \frac{P(y^2)}{\sqrt{U(y^2) V^{3}(y^{2})}} \text{   ,   } g''(y) = \frac{y R(y^2)}{U(y^2) V(y^2) \sqrt{U(y^2) V^{3}(y^2)}},
\ees
where
\be\label{decomp_abcdphase}
\begin{aligned}
&P(z) = 1 - 2(a+c)z + (3ac-bd-(a+c)(b+d))z^2 + 2ac(b+d) z^3 + abcd z^4,\\
&U(z)=(1- a z)(1- c z),\\
&V(z)=(1+ b z)(1+ d z),\\
&R(z) = 2P'(z)U(z)V(z)-P(z)U'(z)V(z)-3P(z)U(z)V'(z).
\end{aligned}
\ee
Note also that $\deg(R) \leq 6$, $\deg(R) \leq 5$ if $bd=0$ and $\deg(R) \leq 4$ if $ac=0$. 

\subsubsection{Low frequencies}

In this case the decay depends on the behavior of $g'$ close to $0$.

\bpr
Let $a,b,c,d$ satisfying \eqref{condition_abcd}-\eqref{condition_abcd2}. There exists a number $y_{0} >0$ such that for any compactly supported function $\chi$ that is equal to $1$ near $0$ and whose support is a subset of $[-y_{0},y_{0}]$, there exists $C>0$ such that for any $\mu \in (0,1]$ and any $t \in \R^{\ast}$
\begin{enumerate}
\item if $a+b+c+d \neq 0$
\begin{align*}
&\left\|e^{\iD \tfrac{t}{\sqrt{\mu}} g(\sqrt{\mu} |D|)} \chi(\sqrt{\mu} D) u_{0} \right\|_{L^{\infty}_{x}} \leq C \frac{\| u_{0}\|_{L^{1}}}{(\mu |t|)^{\frac13}} \text{  ,   when } n=1,\\
&\left\|e^{\iD \tfrac{t}{\sqrt{\mu}} g(\sqrt{\mu} |D|)} \chi(\sqrt{\mu}  |D|) u_{0} \right\|_{L^{\infty}_{x}} \leq C \frac{\| u_{0}\|_{L^{1}}}{\sqrt{\mu} |t|} \text{  ,   when } n=2,
\end{align*}
\item if $a+b+c+d = 0$
\begin{align*}
&\left\|e^{\iD \tfrac{t}{\sqrt{\mu}} g(\sqrt{\mu} |D|)} \chi(\sqrt{\mu} D) u_{0} \right\|_{L^{\infty}_{x}} \leq C \frac{\| u_{0}\|_{L^{1}}}{(\mu^2 |t|)^{\frac15}} \text{  ,   when } n=1,\\
&\left\|e^{\iD \tfrac{t}{\sqrt{\mu}} g(\sqrt{\mu} |D|)} \chi(\sqrt{\mu}  |D|) u_{0} \right\|_{L^{\infty}_{x}} \leq C \frac{\| u_{0}\|_{L^{1}}}{\mu^{\tfrac{3}{5}} |t|^{\tfrac{4}{5}}} \text{  ,   when } n=2,
\end{align*}
\item in any case when $n=2$ for any $k \in \Z$ such that $2^{k} \leq \frac{y_{0}}{2}$
\[
\left\|e^{\iD \tfrac{t}{\sqrt{\mu}} g(\sqrt{\mu} |D|)} Q_{k}(\sqrt{\mu} |D|) u_{0} \right\|_{L^{\infty}_{x}} \leq C \frac{\| |D|^{\frac32}  Q_{k}(\sqrt{\mu} |D|) u_{0} \|_{L^{1}}}{\sqrt{|t|}}.
\]
\end{enumerate}
\epr

Corresponding Strichartz estimates can also be obtained from Propositions \ref{Strichartz_chi} and \ref{Strichartz_Paleysum}.

\begin{proof}
We notice that $g'(y) \underset{0}{=} 1 - \tfrac32 (a+b+c+d) y^2 + O(y^4)$ and if $a+b+c+d=0$, $g'(y) \underset{0}{=} 1 + \tfrac52 (a+b)(b+c) y^4 + O(y^6)$. The bounds follow from Lemma \ref{decay1d_VdC}, Lemma \ref{decay2d_ww_lowfreq} and Lemma \ref{decay2d_wave_2deriv}.
\end{proof}

\subsubsection{Intermediate frequencies}

In that case, the decay depends on how many derivatives of $g$ vanish at the same time so that one has to check the roots of the polynomials $P$ and $R$ defined in \eqref{decomp_abcdphase}. The key observation is that $g', g'', \cdots , g^{(l)}$ have a common zero if and only if $P,P',\cdots,P^{(l-1)}$ have a common root and that $g'', g''', \cdots , g^{(l)}$ have a common zero if and only if $R,R',\cdots,R^{(l-2)}$ have a common root. We begin with the case $n=1$.

\bpr
Let $a,b,c,d$ satisfying \eqref{condition_abcd}-\eqref{condition_abcd2} and $n= 1$. Define $R$ as in \eqref{decomp_abcdphase} and $m \in \N$ as the maximum among the multiplicities of positive roots of $R$. Then for any $y_{1}>y_{0}>0$, any smooth function $\tilde \chi$ that is compactly supported in $[y_{0},y_{1}]$, there exists $C>0$ such that for any $\mu \in (0,1]$ and any $t \in \R^{\ast}$
\[
\left\|e^{\iD \tfrac{t}{\sqrt{\mu}} g(\sqrt{\mu} |D|)} \tilde \chi(\sqrt{\mu}  |D|) u_{0} \right\|_{L^{\infty}_{x}} \leq C \frac{\| u_{0} \|_{L^{1}}}{\mu^{\frac{m+1}{2m+4}} |t|^{\frac{1}{m+2}}}.
\]
\epr
Corresponding Strichartz estimates can also be obtained from Proposition \ref{Strichartz_chi}.

\br\label{comment_n=1_midfreqabcd}
Note that $m \leq \deg(R) \leq 6$. In many cases one can compute $m$, at least numerically. For instance, in the case $a=b=c=0$ and $d>0$, $R(r)=-3d$ and $m=0$.
\er

\begin{proof}
This is a consequence of Lemma \ref{decay1d_multideriv} since $\sum_{k=0}^{m} |R^{(k)}| > 0$ on $\R^{\ast}_{+}$ by definition of $m$.
\end{proof}

Before stating the result when $n=2$, we have the following technical lemma.

\bl 
Let $a,b,c,d$ satisfying \eqref{condition_abcd}-\eqref{condition_abcd2}. Using the notations from \eqref{decomp_abcdphase}, $P$, $P'$ and $P''$ do not have a common root.
\el

\begin{proof}
Arguing by contradiction, there exists $\lambda>0$, $u,v \in \R$ such that $(X-\lambda)^3(u X + v) = P$. By checking the $X^{4}$, $X^{0}$ and then $X^{3}$ coefficients, we get $u=abcd$, $v=-\tfrac{1}{\lambda^{3}}$ and then 
\[
1+2ac(b+d)\lambda^{3} + 3 abcd \lambda^{4} = 0 
\] 
which is impossible since $ac(b+d) \geq 0$ and $abcd \geq 0$.
\end{proof}

One can now state a result in the case $n=2$.
\bpr
Let $a,b,c,d$ satisfying \eqref{condition_abcd}-\eqref{condition_abcd2} and $n= 2$. Define $P,R$ as in \eqref{decomp_abcdphase} and $m \in \N$ as the maximum among the multiplicities of positive roots of $R$. Then, for any $y_{1}>y_{0}>0$, any smooth function $\tilde \chi$ that is compactly supported in $[y_{0},y_{1}]$, there exists $C>0$ such that for any $\mu \in (0,1]$ and any $t \in \R^{\ast}$

\begin{enumerate}
\item if $P$ does not have a positive root
\[
\left\|e^{\iD \tfrac{t}{\sqrt{\mu}} g(\sqrt{\mu} |D|)} \tilde \chi(\sqrt{\mu}  |D|) u_{0} \right\|_{L^{\infty}_{x}}  \leq C \min \left( \frac{\| u_{0} \|_{L^{1}}}{\mu^{\frac{3m+4}{4m+8}} |t|^{\frac{m+4}{2m+4}}} , \frac{\| |D|^{\frac32} \tilde \chi(\sqrt{\mu}  |D|) u_{0} \|_{L^{1}}}{\sqrt{|t|}} \right),
\]
\item if $P$ has a positive root and if $P$ and $P'$ do not have a common positive root
\[
\left\|e^{\iD \tfrac{t}{\sqrt{\mu}} g(\sqrt{\mu} |D|)} \tilde \chi(\sqrt{\mu}  |D|) u_{0} \right\|_{L^{\infty}_{x}}  \leq \frac{C}{\sqrt{|t|}} \min \left( \frac{1}{\mu^{\frac{3}{4}}} \|  u_{0}\|_{L^{1}}, \| |D|^{\frac32}  \tilde \chi(\sqrt{\mu}  |D|) u_{0} \|_{L^{1}} \right),
\]
\item if $P$ and $P'$ have a common positive root
\[
\left\|e^{\iD \tfrac{t}{\sqrt{\mu}} g(\sqrt{\mu} |D|)} \tilde \chi(\sqrt{\mu}  |D|) u_{0} \right\|_{L^{\infty}_{x}}  \leq \frac{C}{|t|^{\frac13}} \min \left( \frac{1}{\mu^{\frac{5}{6}}} \| u_{0} \|_{L^{1}} , \| |D|^{\frac53} \tilde  \chi(\sqrt{\mu}  |D|) u_{0} \|_{L^{1}} \right).
\]
\end{enumerate}

\epr
Corresponding Strichartz estimates can also be obtained from Propositions \ref{Strichartz_chi} and \ref{Strichartz_Paleysum}. 

\begin{proof}
This is a consequence of  Lemma \ref{decay2d_lowfreq_g'not0}, Remark \ref{decay2d_lowfreq_g'cancels}, Lemma \ref{decay2d_wave_multideriv} since $\sum_{k=0}^{m} |R^{(k)}| > 0$ on $\R^{\ast}_{+}$ by definition of $m$ and that $|P|+|P'|+ |P''|>0$ on $\R^{+}$ by the previous lemma. 
\end{proof}

\subsubsection{High frequencies} After careful computations one can show that if $a,b,c,d$ satisfies \eqref{condition_abcd}-\eqref{condition_abcd2} we have
\ba\label{boundshighfreqabcd}
g'(y) - \ell \underset{\infty}{\sim} \Gamma_{1} y^{\alpha+1} \text{  ,  } g''(y) \underset{\infty}{\sim} (\alpha+1) \Gamma_{1} y^{\alpha} \text{  ,  } g'''(y) \leq \Gamma_{2} y^{\alpha-1},
\ea
where $\Gamma_{1}$ and $\Gamma_{2}$ depend on $a,b,c,d$ and $\alpha,\ell$ are given by Table \ref{tablehighfreqabcd}.

We obtain from Lemma \ref{decay1d_VdC2}, Lemma \ref{VdC2d}, Lemma \ref{decay2d_theta=0_HF} and Proposition \ref{Strichartz_Paleysum} the following proposition.
\bpr\label{highfreqBoussi}
Let $a,b,c,d$ satisfying \eqref{condition_abcd}-\eqref{condition_abcd2} and $n=1$ or $2$. Define $\alpha$ and $\ell$ as in Table \ref{tablehighfreqabcd}. There exists a number $y_{1} >0$ such that for any compactly supported function $\chi$ that is equal to $1$ on $[-y_{1},y_{1}]$, any $(q,r)$ sharp $\tfrac{n}{2}$-admissible, there exists $C>0$ such that for any $\mu \in (0,1]$ and any $t \in \R^{\ast}$
\begin{enumerate}
\item if $n=2$ and $\ell \neq 0$
\begin{align*}
&\left\|e^{\iD \tfrac{t}{\sqrt{\mu}} g(\sqrt{\mu} |D|)} (1-\chi(\sqrt{\mu}  |D|)) u_{0} \right\|_{L^{\infty}_{x}} \leq C \mu^{- \frac{\alpha+1}{4}} \frac{\| |D|^{\frac{1-\alpha}{2}} (1-\chi(\sqrt{\mu} |D|)) u_{0} \|_{L^{1}}}{|t|},\\
&\left\|e^{\iD \tfrac{t}{\sqrt{\mu}} g(\sqrt{\mu} |D|)} (1-\chi(\sqrt{\mu}  |D|)) u_{0} \right\|_{L^{q}_{t} L^{r}_{x}} \leq C \mu^{- \frac{\alpha+1}{4}(\frac12-\frac1r)} \| |D|^{\frac{1-\alpha}{2}(\frac12-\frac1r)} (1-\chi(\sqrt{\mu} |D|)) u_{0} \|_{L^{2}},
\end{align*}
\item else
\begin{align*}
&\left\|e^{\iD \tfrac{t}{\sqrt{\mu}} g(\sqrt{\mu} |D|)} (1-\chi(\sqrt{\mu}  |D|)) u_{0} \right\|_{L^{\infty}_{x}} \leq C \mu^{- \frac{n}{2} \frac{\alpha+1}{2}} \frac{\| |D|^{-\frac{n}{2} \alpha} (1-\chi(\sqrt{\mu} |D|)) u_{0} \|_{L^{1}}}{|t|^{\frac{n}{2}}},\\
&\left\|e^{\iD \tfrac{t}{\sqrt{\mu}} g(\sqrt{\mu} |D|)} (1-\chi(\sqrt{\mu}  |D|)) u_{0} \right\|_{L^{q}_{t} L^{r}_{x}} \leq C \mu^{- \frac{n}{2} \frac{\alpha+1}{2} (\frac12-\frac1r)} \| |D|^{-\frac{n}{2} \alpha (\frac12-\frac1r)} (1-\chi(\sqrt{\mu} |D|)) u_{0} \|_{L^{2}}.
\end{align*}
\end{enumerate}
\epr 

\begin{table}[!h]
\centering
\begin{tabular}{|c|c|c|}
  \hline
  & $\ell$ & $\alpha$ \\\hline
  $b=d=0$, $a<0$,  $c<0$ & $0$ & $1$ \\\hline
  $b=d=0$, $ac = 0$, $a+c<0$ & $0$ & $0$ \\\hline
  $bd=0$, $b+d>0$, $ac \neq 0$ & $0$ & $0$ \\\hline
  $bd=0$, $b+d>0$, $ac = 0$, $a+c<0$ & $\sqrt{\frac{-(a+c)}{b+d}}$ & $-3$ \\\hline
  $bd=0$, $b+d>0$, $a=c = 0$ & $0$ & $-4$ \\\hline
  $bd \neq 0$, $a=c = 0$ & $0$ & $-3$ \\\hline
  $bd \neq 0$, $ac = 0$, $a+c<0$, $b(a+c)+bd+(a+c)d \neq 0$ & $0$ & $-4$ \\\hline
  $bd \neq 0$, $ac = 0$, $a+c<0$, $b(a+c)+bd+(a+c)d = 0$ & $0$ & $-6$ \\\hline
  $abcd \neq 0$, $abc+abd+acd+bcd \neq 0$ & $\sqrt{\frac{ac}{bd}}$ & $-3$ \\\hline
  $abcd \neq 0$, $abc+abd+acd+bcd = 0$  & $\sqrt{\frac{ac}{bd}}$ & $-5$ \\\hline
\end{tabular}
\caption{Values of $\alpha$ and $\ell$ in \eqref{boundshighfreqabcd} with $a,b,c,d$ under Assumptions \eqref{condition_abcd}-\eqref{condition_abcd2}.}
\label{tablehighfreqabcd}
\end{table}

\subsubsection{Global Strichartz estimates}

One can gather the previous bounds on the low, intermediate and high frequencies in order to prove global $L^{\infty}$-decays and global Strichartz estimates. For instance, in the case $n=2$, $ac>0$ and $b=d=0$, introducing a smooth compactly supported function $\chi$ that is equal to $1$ near $0$, we obtain from the previous subsection that there exists $C>0$ such that for any $\mu \in (0,1]$ and any $t \in \R^{\ast}$
\begin{align*}
\left\|e^{\iD \tfrac{t}{\sqrt{\mu}} g(\sqrt{\mu} |D|)} u_{0} \right\|_{L^{\infty}_{x}} \leq \frac{C}{\mu^{\frac12} |t|} \| \chi(\sqrt{\mu} |D|) u_{0} \|_{L^{1}} + \frac{C}{\mu |t|} \| |D|^{-1} (1-\chi(\sqrt{\mu} |D|)) u_{0} \|_{L^{1}}
\end{align*}
improving the decay provided by \cite{LPS12} (Proposition 2). In the following we focus on global Strichartz estimates, since $L^{\infty}$ decays are obtained in the same way. We begin with the case $n=1$ which is a consequence of Lemma \ref{decay1d_VdC2}, Lemma \ref{decay1d_multideriv} and Proposition \ref{Strichartz_Paleysum}.

\bpr 
Let $a,b,c,d$ satisfying \eqref{condition_abcd}-\eqref{condition_abcd2} and $n= 1$. Define $R$ as in \eqref{decomp_abcdphase} and $m \in \N$ as the maximum among the multiplicities of positive roots of $R$. Let $l=\max(m+2,3)$ if $a+b+c+d \neq 0$ or $l=\max(m+2,5)$ if $a+b+c+d = 0$. Define $\alpha$ as in Table \ref{tablehighfreqabcd} and a smooth compactly supported function $\chi$ that is equal to $1$ near $0$. Then for any $(q,r) \in [2,\infty]^2$ with $(q,r)$ sharp $\frac{1}{l}$-admissible, there exists a constant $C>0$ such that for any $\mu \in (0,1]$
\begin{align*}
\left\|e^{\iD \tfrac{t}{\sqrt{\mu}} g(\sqrt{\mu} |D|)} u_{0} \right\|_{L^{q}_{t} L^{r}_{x}} \leq &C  \mu^{\frac{1-l}{2l}(\frac12-\frac1r)} \| \chi(\sqrt{\mu} |D|) u_{0} \|_{L^{2}}\\
& + C \mu^{-\frac{1+\alpha}{2l} (\frac12-\frac1r)} \| |D|^{\frac{l-2-\alpha}{l} (\frac12-\frac1r)} (1-\chi(\sqrt{\mu} |D|)) u_{0} \|_{L^{2}}.
\end{align*}
\epr 
Retarded Strichartz estimates can also be obtained. We now deal with the case $n=2$ which is a consequence of almost all the lemmas in Subsection \ref{ss:n=2} and propositions in Subsection \ref{ss:Strichartz}.

\bpr 
Let $a,b,c,d$ satisfying \eqref{condition_abcd}-\eqref{condition_abcd2} and $n= 2$. Define $P,R$ as in \eqref{decomp_abcdphase} and $m \in \N$ as the maximum among the multiplicities of positive roots of $R$. Let $u = \frac{2m+4}{m+4}$ if $P>0$ on $\R^{+}$ or $u=2$ if $P$ has a positive root and if $|P|+|P'|>0$ on $\R^{+}$ or $u=3$ if $P$ and $P'$ have a common positive root. Let $k=1$ if $a+b+c+d \neq 0$ or $k=\frac{5}{4}$ if $a+b+c+d = 0$. Let $l=\max(u,k)$. Define $\alpha,\ell$ as in Table \ref{tablehighfreqabcd} and a smooth compactly supported function $\chi$ that is equal to $1$ near $0$. Then for any $(q,r) \in [2,\infty]^2$ with $(q,r)$ sharp $\frac{1}{l}$-admissible, there exists $C>0$ such that for any $\mu \in (0,1]$
\begin{enumerate}
\item if $\ell \neq 0$ and additionally when $l>1$, $\alpha \neq -3$
\begin{align*}
\left\|e^{\iD \tfrac{t}{\sqrt{\mu}} g(\sqrt{\mu} |D|)} u_{0} \right\|_{L^{q}_{t} L^{r}_{x}} \leq &C  \mu^{\frac{1-2l}{2l}(\frac12-\frac1r)} \| \chi(\sqrt{\mu} |D|) u_{0} \|_{L^{2}}\\
& + C \mu^{-\frac{1+\alpha}{4l} (\frac12-\frac1r)} \| |D|^{\frac{4l-3-\alpha}{2l} (\frac12-\frac1r)} (1-\chi(\sqrt{\mu} |D|)) u_{0} \|_{L^{2}},
\end{align*}
\item if $\ell \neq 0$, $\alpha=-3$ and $l>1$, we  must assume $r<\infty$ and
\begin{align*}
\left\|e^{\iD \tfrac{t}{\sqrt{\mu}} g(\sqrt{\mu} |D|)} u_{0} \right\|_{L^{q}_{t} L^{r}_{x}} \leq &C  \mu^{\frac{1-2l}{2l}(\frac12-\frac1r)} \| \chi(\sqrt{\mu} |D|) u_{0} \|_{L^{2}}\\
& + C \mu^{-\frac{l-1}{4l} (\frac12-\frac1r)} \| |D|^{(2-\frac{1}{l}) (\frac12-\frac1r)} (1-\chi(\sqrt{\mu} |D|)) u_{0} \|_{L^{2}},
\end{align*}
\item else
\begin{align*}
\left\|e^{\iD \tfrac{t}{\sqrt{\mu}} g(\sqrt{\mu} |D|)} u_{0} \right\|_{L^{q}_{t} L^{r}_{x}} \leq &C  \mu^{\frac{1-2l}{2l}(\frac12-\frac1r)} \| \chi(\sqrt{\mu} |D|) u_{0} \|_{L^{2}}\\
& + C \mu^{-\frac{1+\alpha}{2l} (\frac12-\frac1r)} \| |D|^{\frac{2l-2-\alpha}{l} (\frac12-\frac1r)} (1-\chi(\sqrt{\mu} |D|)) u_{0} \|_{L^{2}}.
\end{align*}
\end{enumerate}
\epr 
Retarded Strichartz estimates can also be obtained as well as Strichartz estimates that are uniform with respect to $\mu \in (0,1]$.

\subsubsection{Local Kato smoothing/Morawetz type estimates}

Let $a,b,c,d$ satisfying \eqref{condition_abcd}-\eqref{condition_abcd2}, $n =1$ or $2$ and define $\alpha$ and $\ell$ as in Table \ref{tablehighfreqabcd}.

If $\alpha \geq 0$, there is a local Kato smoothing. Let $\chi$ be a compactly supported function as in Proposition \ref{highfreqBoussi}. Using Proposition \ref{kato-morawetz} there exists $C>0$ such that for any $\mu \in (0,1]$
\[
\int_{\R} \int_{\R^{n} } e^{-x^{2}} |(e^{\iD \frac{t}{\sqrt{\mu}} g(\sqrt{\mu} |D|)} |D|^{\frac{\alpha+1}{2}} (1-\chi(\sqrt{\mu}  |D|)) u_{0})(x)|^2 dx dt \leq C \mu^{- \frac{\alpha+1}{2}} \left\| u_{0} \right\|_{L^{2}}^2.
\]

If $\alpha < 0$ and $P>0$, we can prove Morawetz type estimates. From Proposition \ref{kato-morawetz} there exists $C>0$ such that for any $\mu \in (0,1]$
\begin{align*}
&\sup_{x_{0} \in \R^{n}} \int_{\R} \int_{\R^{n} } |(e^{\iD \frac{t}{\sqrt{\mu}} g(\sqrt{\mu} |D|)} u_{0})(x)|^2 e^{-|x-x_{0}|^{2}} dx dt \leq C \left\| u_{0} \right\|_{L^{2}}^2 \text{   ,   in case } \ell \neq 0,\\
&\sup_{x_{0} \in \R^{n}} \int_{\R} \int_{\R^{n} } |(e^{\iD \frac{t}{\sqrt{\mu}} g(\sqrt{\mu} |D|)} u_{0})(x)|^2 e^{-|x-x_{0}|^{2}} dx dt \leq C \left\| (1+\sqrt{\mu} |D|)^{-\frac{\alpha+1}{2}} u_{0} \right\|_{L^{2}}^2 \text{   ,   in case } \ell = 0.
\end{align*}
Corresponding decay of local energy can be obtained from Corollary \ref{decay_local_energy}.

\bibliographystyle{alpha}
\bibliography{biblio}
\end{document}